\documentclass[12pt,a4paper]{amsart}
\usepackage{amssymb,amscd,amsmath,mathrsfs}
\usepackage{times}

\input xy
\xyoption{all}

\newif\ifEU\IfFileExists{euscript.sty}{\EUtrue\usepackage[mathcal]{euscript}}{}

\makeatletter

\renewcommand{\phi}{\varphi}
\renewcommand{\to}{\longrightarrow}
\renewcommand{\le}{\leqslant}
\renewcommand{\geq}{\geqslant}
\renewcommand{\epsilon}{\varepsilon}

\renewcommand{\kappa}{\varkappa}
\DeclareMathOperator{\CAT}{\mathbf{CAT}}
\DeclareMathOperator{\Cyl}{Cyl}
\DeclareMathOperator{\Cocyl}{Cocyl}

\DeclareMathOperator{\Sets}{Sets}

 \DeclareMathOperator{\Hom}{Hom}

\DeclareMathOperator{\id}{id} \DeclareMathOperator{\Cat}{{\bf
Cat}}

 \DeclareMathOperator{\Ar}{Ar}

\DeclareMathOperator{\Ob}{Ob}
\newcommand{\la}{\textrm{\huge${}_{\ulcorner}$}}

\newcommand{\lra}[1]{\bl{#1}\longrightarrow\relax}
\newcommand{\bl}[1]{\buildrel #1\over}
\newcommand{\cc}{\mathscr}
\newcommand{\bb}{\mathbf}
\newcommand{\ps}{\oplus}

\newcommand{\iso}{\simeq}
\newcommand{\op}{{\textrm{\rm op}}}

\newcommand{\wt}{\widetilde}

\newcommand {\di}{\mathcal Dia}
\newcommand {\ord}{\mathcal Ord}

\newcommand {\dirf}{\mathcal{D}\mathit{irf}}
\newcommand{\homa}[2]{\Hom(#1,#2)}

\newtheorem{thm}{Theorem}[section]
\newtheorem{prop}[thm]{Proposition}
\newtheorem*{funcax}{\bf\em Functoriality Axiom}
\newtheorem*{isomax}{\bf\em Isomorphism Axiom}
\newtheorem*{disjax}{\bf\em Disjoint Union Axiom}
\newtheorem*{kanax}{\bf\em Homotopy Kan Extension Axiom}
\newtheorem*{baseax}{\bf\em Base Change Axiom}
\newtheorem{cor}[thm]{Corollary}
\newtheorem{lem}[thm]{Lemma}

\newtheorem*{lemb}{Lemma B}

\newtheorem*{question}{Question}
\newtheorem*{sublem}{Sublemma}

\newtheorem*{defs}{Definition}

\makeatother

\begin{document}

\footskip30pt


\title{Systems of Diagram Categories and K-theory. II}
\author{Grigory Garkusha}
\address{International Centre for Theoretical Physics, Strada Costiera 11, I-34014, Trieste, Italy}
\urladdr{www.ictp.trieste.it/$\sim$garkusha}
\email{ggarkusha@mail.ru}
\thanks{Supported by the ICTP Research Fellowship}
\begin{abstract}
The additivity theorem for d\'erivateurs associated to complicial
biWaldhau\-sen categories is proved. Also, to any exact category
in the sense of Quillen a $K$-theory space is associated. This
$K$-theory is shown to satisfy the additivity, approximation and
resolution theorems.
\end{abstract}
\maketitle

\thispagestyle{empty} \pagestyle{plain}

\newdir{ >}{{}*!/-6pt/@{>}} 

\section{Introduction}

It is well known due most recently to work of
Schlichting~\cite{Sch} that, in general, there is no $K$-theory
for triangulated categories satisfying localization and
reconstructing Quillen's $K$-theory of an exact category from its
derived category. There are two approaches to replace the naive
notion of derived category by something richer, from which the
$K$-theory might be obtained by some explicit construction.  One
approach, suggested by Dwyer and Kan~\cite{DK1,DK2}, is called the
``simplicial localization". To\"en and Vezzosi~\cite{TV} define a
$K$-theory functor on the level of $S$-categories and prove that,
when applied to the simplicial localization of an appropriate
Waldhausen category $\cc C$, this construction yields a spectrum
which is weakly equivalent to the Waldhausen's $K$-theory spectrum
of $\cc C$. The other, following definitions of
Grothendieck~\cite{G}, Heller~\cite{H} and Franke~\cite{F} is
called the ``system of diagram categories" or the ``d\'erivateur"
(in French): given a closed model category $\cc C$, one takes a
large system of index categories $I$, and forms the system of
derived categories of the diagram categories $\cc C^I$. One can
introduce an analogue of the $Q.$-construction (see~\cite{M}) or
Waldhausen's $S.$-construction~\cite{Gar} for them which might
form a suitable definition for the $K$-theory of a system of
diagram categories or a d\'erivateur. Both definitions give
equivalent $K$-theories by~\cite{C2}.

Maltsiniotis~\cite{M} opens three conjectures, the third of which
says that the $K$-theory of a d\'erivateur (or a system of diagram
categories) satisfies the additivity theorem. A weaker version of
additivity has been shown in~\cite{Gar}. One could try to adapt
Waldhausen's~\cite{W} or McCarthy's~\cite{Mc} proof of additivity
in this context to show the strong form of additivity. At the
first glance, the machinery of d\'erivateurs seems to have some of
the necessary ingredients. However one thing always goes wrong
when constructing a homotopy; one lacks a component which seems to
be not deducible from d\'erivateurs in general (see details at the
end of section~\ref{c}).

In this paper we prove the strong form of additivity for
d\'erivateurs associated to complicial biWaldhausen categories in
the sense of Thomason~\cite{T} (see the precise definitions in
section~\ref{c}). Thus we can find the lacking ingredient in this
case. Experts will probably be able to show additivity for the
d\'erivateurs represented by closed model categories. Such
d\'erivateurs seem to contain all the necessary information for
this. If so, it will be justifiable to say that the third
Maltsiniotis conjecture is true for the d\'erivatuers ``having a
model". In fact, all interesting d\'erivatuers we have in nature
arise in this fashion.

To any exact category $\cc E$ one associates the d\'erivateur
   $$\bb D^b(\cc E):I\longmapsto D^b(\cc E^I)$$
sending an index category $I$ to the derived category of the exact
diagram category $\cc E^I$. It is very interesting to study its
$K$-theory space $K(\bb D^b(\cc E))$. This $K$-theory is shown to
satisfy the additivity, approximation and resolution theorems.

\subsubsection*{Organization of the paper} After fixing some
notation and terminology (in~\ref{2.1} and~\ref{2.2}), we
formulate several lemmas (in~\ref{rid}) which are of great utility
in proving the additivity theorem in section~\ref{234}. Our main
result is then formulated in section~\ref{c}, that dealing with
additivity for d\'erivateurs associated to complicial
biWaldhau\-sen categories (Theorem~\ref{mod}). Then comes
section~\ref{5} in which the $K$-theory space $K(\bb D^b(\cc E))$
is studied. We also prove there a couple of results which are of
independent interest. The necessary facts about d\'erivateurs and
their $K$-theory are given in Addendum.

\subsubsection*{Acknowledgement} I would like to thank Professor
Haynes Miller and an anonymous referee for helpful suggestions
concerning the material of the paper.

\section{Cochain complexes, homotopies, derived categories}

In this section we collect some necessary facts about cochain
complexes and derived categories.

\subsection{Definition of the derived category}\label{2.1}

An {\it exact category\/} is an additive category $\cc A$ with a
collection of {\it exact sequences\/} $\{E\rightarrowtail
F\twoheadrightarrow G\}$ where the first morphism
$E\rightarrowtail F$ appearing in those exact sequences are called
{\it admissible monomorphisms\/} and the second ones {\it
admissible epimorphisms}. They have to satisfy a couple of natural
axioms (e.g. see~\cite{T}). Let $\cc C=C^b(\cc A)$ denote the
category of bounded complexes in an exact category $\cc A$. Recall
that $\cc A$ can be embedded as a full subcategory in an abelian
category $\cc B$ in such a way that a sequence in $\cc A$ is exact
if and only if it is exact in $\cc B$. If $\cc A$ is idempotent
complete (or even less) this embedding can be chosen in a way that
any map in $\cc A$ which becomes an epimorphism in $\cc B$ was
already an admissible epimorphism in $\cc A$ (see~\cite{T}).

Let $\cc A$ be an exact category. Its bounded derived category
$D^b(\cc A)$ is constructed as follows (we follow here Keller's
definition~\cite{K}).

Let $H^b(\cc A)$ be the homotopy category of the category of
bounded complexes $\cc C=C^b(\cc A)$, i.e., the quotient category
of $\cc C$ modulo homotopy equivalence. Let $Ac(\cc A)$ denote the
full subcategory of $H^b(\cc A)$ consisting of acyclic complexes.
A complex
   $$X^n\lra{}X^{n+1}\lra{}X^{n+2}$$
is called {\it acyclic\/} if each map $X^n\to X^{n+1}$ decomposes
in $\cc A$ as $X^n\bl{e_n}\twoheadrightarrow
D^n\bl{m_n}\rightarrowtail X^{n+1}$ where $e_n$ is an epimorphism
and $m_n$ is a monomorphism in such a way that
$D^n\bl{m_n}\rightarrowtail X^{n+1}\bl{e_{n+1}}\twoheadrightarrow
D^{n+1}$ is an exact sequence.

If an exact category is idempotent complete then every
contractible complex is acyclic. Denote by $\cc N=\cc N(\cc A)$
the full subcategory of $H^b(\cc A)$ whose objects are the
complexes isomorphic in $H^b(\cc A)$ to acyclic complexes. There
is another description of $\cc N$. Let $\cc A\to\tilde{\cc A}$ be
the universal additive functor to an idempotent complete exact
category $\tilde{\cc A}$. It is exact and reflects exact
sequences, and $\cc A$ is closed under extensions in $\tilde{\cc
A}$ (see~\cite[A.9.1]{T}). The class of acyclic complexes in
$\tilde{\cc A}$ is closed under homotopy equivalence. It follows
that a complex with entries in $\cc A$ belongs to $\cc N$ if and
only if its image in $H^b(\tilde{\cc A})$ is acyclic. The category
$\cc N(\tilde{\cc A})=Ac(\tilde{\cc A})$ is a thick subcategory in
$H^b(\tilde{\cc A})$. Note that a complex over $\tilde{\cc A}$ is
acyclic if and only if it has trivial homology computed in an
appropriate ambient abelian category $\cc B$ (see above). It
follows that $\cc N$ is a thick subcategory in $H^b(\cc A)$.
Denote by $\Sigma$ the multiplicative system associated to $\cc N$
and call the elements of $\Sigma$ {\it quasi-isomorphisms}. A map
$s$ is a quasi-isomorphism if and only if in any triangle
   $$L\lra{s}M\to N\to L[1]$$
the complex $N$ belongs to $\cc N$.

The derived category is defined as
   $$D^b(\cc A)=H^b(\cc A)/\cc N=H^b(\cc A)[\Sigma^{-1}].$$
Clearly, a map is a quasi-isomorphism if and only if its image in
$C^b(\tilde{\cc A})$ is a quasi-isomorphism and if and only if its
image in $D^b(\cc A)$ is an isomorphism.

We shall work a lot with derived categories of diagram exact
categories $\cc A^I$ where $I$ is a small category. It is easily
seen that a cochain map $f:A\to B$ in $C^b(\cc A^I)$ is a
quasi-isomorphisms if and only if each $f_i:A_i\to B_i, i\in I$,
is so in $C^b(\cc A)$.

\subsection{Homotopy pullbacks and homotopy pushouts}\label{2.2}

Let $f:F\to A$ and $g:G\to A$ be cochain maps. One has a canonical
homotopy pullback
   \begin{gather*}
    \Bigl(F\prod_A^h G\Bigr)^n=F^n\oplus A^{n-1}\oplus G^n\\
    d(x,a,y)=(d_Fx,-d_Aa+fx-gy,d_Gy).
   \end{gather*}
(We describe $d$ as if objects of $\cc A$ had ``elements", by the
standard abuse). The square
   $$\xymatrix{\ar @{}[dr] |{\textrm{hom. comm.}}
               F\prod\limits_A^h G \ar[d]_{g'} \ar[r]^{f'} & G \ar[d]^{g} \\
               F \ar[r]_{f}        & A}$$
with $f'(x,a,y)=y$ and $g'(x,a,y)=x$ is homotopy commutative, that
is $gf'\sim fg'$. Note that $f'$ is a quasi-isomorphism whenever
$f$ is. Cochain maps from a complex $C$ to this canonically
homotopy pullback correspond bijectively to data $(h,p,k)$ where
$h:C\to F$ and $p:C\to G$ are cochain maps and $k$ is a cochain
homotopy $fh\sim gp:C\to A$. Thus $k$ consists of maps $C^n\to
A^{n-1}$ for all $n$ such that $dk+kd=fh-gp$. To $(h,p,k)$
corresponds the cochain map $t:C\to F\prod\limits_A^h G$ defined
as $t(c)=(hc,kc,pc)$. Then $f't=p$ and $g't=h$. When $f:F\to A$ is
the identity map, the canonically homotopy pullback is the mapping
cocylinder $\Cocyl(g)$ of $g:G\to A$.

Dually, given $f:A\to F$ and $g:A\to G$ the canonically homotopy
pushout is the complex defined by
   \begin{gather*}
    \Bigl(F\coprod_A^h G\Bigr)^n=F^n\oplus A^{n+1}\oplus G^n\\
    d(x,a,y)=(d_Fx+fa,-d_Aa,d_Gy-ga).
   \end{gather*}
This indeed has all the dual properties as the homotopy pullback.
As special cases, when $f:A\to F$ is the identity map, the
homotopy pushout is the mapping cylinder $\Cyl(g)$ of $g:A\to G$.
If $f:A\to F=0$ is the map to 0, the homotopy pushout is the
mapping cone $C(g)$ of $g:A\to G$.

Let $w\cc C$ denote the category whose objects are those of $\cc
C$ and morphisms are quasi-isomorphisms. It is a complical
biWaldhausen category (see definitions in~\cite{T,W}). It has also
cylinder and cocylinder functors satisfying the cylinder and
cocylinder axioms.

\subsection{Getting rid of homotopy commutative squares}\label{rid}

Results of this technical paragraph are of great utility in
proving the ``additivity theorem" in the next section. To
construct a homotopy in that proof we will want to replace some
homotopy commutative diagrams by strictly commutative ones. Given
a non-negative integer $n$, by $\Delta^n$ denote the totally
ordered set $\{0<1<\cdots<n\}$.

Suppose we are given a homotopy commutative square with entries
$(X_0,Y,A_0,A_1)$
   $$\xymatrix{&& X_1\ar@{.>}[dl]^l\ar@/^/@{.>}[ddl]^{g''}\\
               \ar@{}[dr] |{\textrm{h. comm.}}
               X_0 \ar@/^/@{.>}[urr]^{f''}\ar[d]_{g'} \ar[r]^{f'} & Y \ar[d]^{g}\\
               A_0 \ar[r]_{f} & A_1}$$
We want to replace it by a strictly commutative square with
entries $(X_0,X_1,A_0,A_1)$.

Let $X_1=\Cocyl(g)$; then $gl\sim g''$. Note that a homotopy is
given by the maps $z^n:X_1^n\to A_1^{n-1}$ mapping $(x,a,y)\in
X_1^n$ to $a$. Since $gf'\sim fg'$ there is a map $f'':X_0\to X_1$
such that $lf''=f'$ and $g''f''=fg'$. A cochain map $f''$ is
defined by $f''(x)=(f'(x),k(x),fg'(x))$ where maps $k^n:X_0^n\to
A_1^{n-1}$ give a homotopy $gf'\sim fg'$.

\begin{lem}\label{mz}
Suppose that in the diagram
   $$\xymatrix{X_0\ar@<2.5pt>[d]_{h'\ }\ar@<-2pt>[d]^{\ \,g'}\ar[r]^{f'}&Y\ar[d]^g\\
               A_0\ar[r]^f&A_1}$$
the maps $h',g'$ are homotopic and the square with $g'$ deleted is
genuinely commutative. Then one can produce a pair of genuinely
commutative squares
   $$\xymatrix{X_0\ar@<2.5pt>[d]_{h'\ }\ar@<-2pt>[d]^{\ \,g'}\ar[r]^{f''}
               &X_1\ar@<2.5pt>[d]_{gl\ }\ar@<-2pt>[d]^{\ \,g''}\\
               A_0\ar[r]^f&A_1}$$
such that the map $(g',g''):X=(X_0\lra{f''}X_1)\to
A=(A_0\lra{f}A_1)$ is homotopic to the map $(h',gl):X\to A$ in
$\cc C^{\Delta^1}=C^b(\cc A^{\Delta^1})$.
\end{lem}

\begin{proof}
Since $fg'\sim gf'$ one can construct a diagram as above with
$X_1=\Cocyl(g)$. By construction, $f''(x)=(f'(x),fm(x),fg'(x))$
where maps $m^n:X_0^n\to A_0^{n-1}$ yield a homotopy $g'\sim h'$.
For any $n$ the square
   $$\begin{CD}
      X_0^n@>{f''^n}>>X_1^n\\
      @V{m^n}VV@VV{z^n}V\\
      A_0^{n-1}@>{f^{n-1}}>>A_1^{n-1}
     \end{CD}$$
is commutative and the maps $(m^n,z^n):X^n\to A^{n-1}$ give the
desired homotopy.
\end{proof}

A map $X\to Y$ in $D^b(\cc A)$ is the equivalence class of a
diagram in $C^b(\cc A)$
   $$X\bl s\longleftarrow Z\lra{f} Y$$
with $s$ a quasi-isomorphism. It is equivalent to $X\bl
t\longleftarrow W\lra{g} Y$ if these fit into a homotopy
commutative diagram
   $$\xymatrix{&&V\ar[dl]_u\ar[dr]^v\\
      &Z\ar[dl]_s\ar[drrr]^(.2)f&& W\ar[dlll]_(.2)t\ar[dr]^g\\
      X&&&&Y}$$
with $u$ and $v$ quasi-isomorphisms.

\begin{lem}\label{imp}
Let $fs^{-1}:A\to C$ be a map in $D^b(\cc A^{\Delta^1})$
represented by a commutative diagram
   \begin{equation}\label{11}
    \begin{CD}
      A_1@<{s_1}<<Y_1@>{f_1}>>C_1\\
      @A{a}AA@AA{y}A@AA{c}A\\
      A_0@<{s_0}<<Y_0@>{f_0}>>C_0
    \end{CD}
   \end{equation}
with $s_0,s_1$ quasi-isomorphisms and let $A_0\bl{t}\longleftarrow
U\lra{h}C_0$ be another representative for
$A_0\xrightarrow{f_0s_0^{-1}}C_0$ in $D^b(\cc A)$ with a common
denominator
   $$
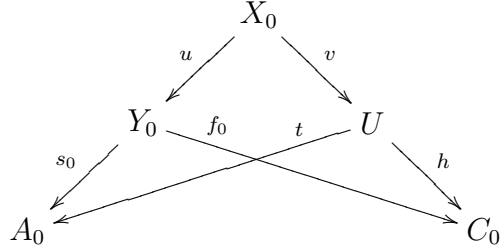
\xymatrix{&&X_0\ar[dl]_u\ar[dr]^v\\
      &Y_0\ar[dl]_{s_0}\ar[drrr]^(.2){f_0}&& U\ar[dlll]_(.2)t\ar[dr]^h\\
      A_0&&&&C_0}$$
Then there exists a complex $X_1$ and a commutative diagram
   \begin{equation}\label{22}
     \begin{CD}
      A_1@<{q}<<X_1@>{g}>>C_1\\
      @A{a}AA@AA{x}A@AA{c}A\\
      A_0@<{tv}<<X_0@>{hv}>>C_0
     \end{CD}
   \end{equation}
representing the same morphism $fs^{-1}$ in $D^b(\cc
A^{\Delta^1})$. If $f_0,f_1$ are quasi-iso\-mor\-phisms then so is
$g$. Moreover, $X_1$ can be chosen in such a way that $x$ is a
monomorphism in $\cc C$.
\end{lem}

\begin{proof}
Applying the preceding lemma first to the diagram
   $$\xymatrix{X_0\ar@<2.5pt>[d]_{s_0u\ }\ar@<-2pt>[d]^{\ \,tv}\ar[r]^{yu}
               &Y_1\ar[d]^{s_1}\\
               A_0\ar[r]^a&A_1}$$
one obtains a diagram
   $$\xymatrix{X_0\ar@<2.5pt>[d]_{s_0u\ }\ar@<-2pt>[d]^{\ \,tv}\ar[r]^{p_1}
               &Y_2\ar@<2.5pt>[d]_{s_1l_1\ }\ar@<-2pt>[d]^{\ \,q_1}\\
               A_0\ar[r]^a&A_1}$$
and then to the diagram
   $$\xymatrix{X_0\ar@<2.5pt>[d]_{f_0u\ }\ar@<-2pt>[d]^{\ \,hv}\ar[r]^{p_1}
               &Y_2\ar[d]^{f_1l_1}\\
               C_0\ar[r]^c&C_1}$$
resulting a diagram
   $$\xymatrix{X_0\ar@<2.5pt>[d]_{f_0u\ }\ar@<-2pt>[d]^{\ \,hv}\ar[r]^{x}
               &X_1\ar@<2.5pt>[d]_{f_1l_1l_2\ }\ar@<-2pt>[d]^{\ \,g}\\
               C_0\ar[r]^c&C_1}$$
We put $q=q_1l_2$. The diagram~\eqref{22} is constructed. It is
equivalent to~\eqref{11}, hence represents the map $fs^{-1}$ in
$D^b(\cc A^{\Delta^1})$. The fact that the map $g$ is a
quasi-isomorphism if $f_0,f_1$ are is obvious.

Finally, to show that $X_1$ can be chosen in such a way that $x$
is a monomorphism in $\cc C$ it is enough to observe that any
cochain map $z:X_0\to Y_3$ is the composite
$X_0\lra{x}X_1\lra{w}Y_3$ of a monomorphism $x$ followed by a
quasi-isomorphism $w$ and $X_1=\Cyl(z)$.
\end{proof}

Let $\Box$ be the poset $\Delta^1\times\Delta^1$ and let
$\la\subset\Box$ be the subposet $\Box\setminus(1,1)$. Then the
exact diagram category ${\cc A}{\ulcorner}$ consists of the
diagrams in $\cc A$
   $$A_{(1,0)}\longleftarrow A_{(0,0)}\lra{}A_{(0,1)}.$$
Let $\wt{\cc A}\ulcorner$ be the full subcategory in ${\cc
A}\ulcorner$ with $A_{(0,0)}\lra{}A_{(0,1)}$ an admissible
monomorphism. In turn, the exact diagram category ${\cc A}^\Box$
consists of the commutative squares
   $$\xymatrix{A_{(0,0)} \ar[d] \ar[r] & A_{(0,1)} \ar[d] \\
           A_{(1,0)} \ar[r]        & A_{(1,1)}}$$
Denote by $\wt{\cc A}^{\Box}$ the full subcategory in ${\cc
A}^\Box$ with $A_{(0,0)}\lra{}A_{(0,1)}$ and
$A_{(1,0)}\lra{}A_{(1,1)}$ admissible monomorphisms in $\cc A$ and
the square above is cocartesian.

It follows that $\wt{\cc A}{\ulcorner}$ is an exact subcategory of
${\cc A}{\ulcorner}$ and $\wt{\cc A}^{\Box}$ is an exact
subcategory of ${\cc A}^{\Box}$. Therefore one can consider their
derived categories $D^b(\wt{\cc A}^{\Box})$ and $D^b(\wt{\cc
A}{\ulcorner})$. We claim that they are naturally equivalent. To
see this, consider the functor $i_{{}_\ulcorner}^*:\wt{\cc
A}^{\Box}\to\wt{\cc A}{\ulcorner}$ taking a square
   $$\xymatrix{A_{(0,0)}\ar[d] \ar@{ >->}[r] & A_{(0,1)} \ar[d] \\
           A_{(1,0)}\ar@{ >->}[r]        & A_{(1,1)}}$$
to $A_{(1,0)}\leftarrow A_{(0,0)}\rightarrowtail A_{(0,1)}$ as
well as the functor $j:\wt{\cc A}{\ulcorner}\to\wt{\cc A}^{\Box}$
taking a diagram $A_{(1,0)}\leftarrow A_{(0,0)}\rightarrowtail
A_{(0,1)}$ to
   $$\xymatrix{A_{(0,0)}\ar[d] \ar@{ >->}[r] & A_{(0,1)} \ar[d] \\
           A_{(1,0)} \ar@{ >->}[r]        & A_{(1,0)}\coprod_{A_{(0,0)}}A_{(0,1)}}$$
Then $i_{{}_\ulcorner}^*$ and $j$ are exact functors and plainly
mutual inverses with $i_{{}_\ulcorner}^*j=\id$. These induce the
desired equivalence of derived categories.

\begin{cor}\label{van}
Given two squares of cochain complexes
   $$\xymatrix{\ar @{}[dr] |{A}
              A_{(0,0)}\ar[d] \ar@{ >->}[r] & A_{(0,1)} \ar[d] \\
              A_{(1,0)}\ar@{ >->}[r]        & A_{(1,1)}}\ \ \ \ \
   \xymatrix{\ar @{}[dr] |{B}
              B_{(0,0)} \ar[d] \ar@{ >->}[r] & B_{(0,1)} \ar[d] \\
              B_{(1,0)} \ar@{ >->}[r]        & B_{(1,1)}}$$
in $D^b(\wt{\cc A}^{\Box})$ and a morphism
$\alpha:i_{{}_\ulcorner}^*(A)\to i_{{}_\ulcorner}^*(B)$ in
$D^b(\wt{\cc A}{\ulcorner})$, there exists a unique map $a:A\to B$
such that $i_{{}_\ulcorner}^*(a)=\alpha$. If $\alpha$ is an
isomorphism then so is $a$.
\end{cor}

\section{The additivity theorem}\label{234}

In this section we prove a sort of the additivity theorem. It
assumes the role of a basic result in algebraic $K$-theory. We
refer the reader to Staffeldt's work~\cite{St}. The author knows
two proofs of that theorem for Waldhausen's categories: by
Waldhausen~\cite{W} and by~McCarthy~\cite{Mc}. We shall follow
Waldhausen's proof.

Waldhausen~\cite{W1,W} constructs a simplicial exact category
$S.\cc A=\{S_n\cc A\}_{n\geq 0}$ in which the face and the
degeneracy maps are exact functors. Let $\Ar\Delta^n$ be the poset
of pairs $(i,j)$, $0\le i\le j\le n$, where $(i,j)\le(i',j')$ if
and only if $i\le i'$ and $j\le j'$. An object of $S_n\cc A$ is a
functor $A:\Ar\Delta^n\to\cc A$ such that $A_{ii}=0$ and
   $$A_{ij}\to A_{ik}\to A_{jk}$$
is a short exact sequence in $\cc A$ for any $0\le i\le j\le k\le
n$. Observe (exercise!) that a cochain map $f:A\to A'$ in
$C^b(S_n\cc A)$ is a quasi-isomorphism if and only if each
$f_{ij}:A_{ij}\to A'_{ij}$ is so in $C^b(\cc A)$.

For any $n\geq 1$ the exact category $S_n\cc A$ is equivalent to
the exact category $F_{n-1}\cc A$ of composable monomorphisms in
$\cc A$
   $$A_0\rightarrowtail A_1\rightarrowtail\cdots\rightarrowtail A_{n-1}.$$
This equivalence is given by the the exact functor forgetting
quotients.

Denote by $i.\bb S.\cc A$ the bisimplicial set
   $$\Delta^m\times\Delta^n\longmapsto i_m\bb S_n\cc A=i_mD^b(S_n\cc A).$$
The $(m,n)$-simplices are represented by the strings of
isomorphisms in $\bb S_n\cc A=D^b(S_n\cc A)$
   $$A_0\lra{\sim}A_1\lra{\sim}\cdots\lra{\sim}A_m.$$
Note that every exact functor $f:\cc A\to\cc A'$ induces a
simplicial map $f_*:i.\bb S.\cc A\to i.\bb S.\cc A'$. We also
observe that coproduct gives a unitial and associative $H$-space
structure to $|i.\bb S.\cc A|$ via the map
   $$|i.\bb S.\cc A|\times|i.\bb S.\cc A|\lra{\sim}|i.\bb S.\cc A\times i.\bb S.\cc A|\lra{\coprod}|i.\bb S.\cc A|.$$

The category $\bb S_0\cc A$ is the trivial category with one
object and one morphism. Hence the geometric realization $|i.\bb
S_0\cc A|$ is the one-point space. The category $\bb S_1\cc A$ is
isomorphic to the derived category $D^b(\cc A)$. Hence the
category of isomorphisms $i\bb S_1\cc A$ may be identified to
$iD^b(\cc A)$.

Consider $|i.\bb S.\cc A|$. The ``1-skeleton" in the
$S.$-direction is obtained from the ``0-skeleton" (which is
$|i.\bb S_0\cc A|$) by attaching of $|i.\bb S_1\cc
A|\times|\Delta[1]|$ (where $|\Delta[1]|$ denotes the topological
space 1-simplex). It follows that the ``1-skeleton" is naturally
isomorphic to the suspension $S^1\wedge|i.D^b(\cc A)|$. One
obtains an inclusion $S^1\wedge|i.D^b(\cc A)|\to|i.\bb S.\cc A|$,
and by adjointness an inclusion of $|i.D^b(\cc A)|$ into the loop
space of $|i.\bb S.\cc A|$,
   $$|i.D^b(\cc A)|\to\Omega|i.\bb S.\cc A|.$$

We can apply the $S.$-construction to produce a bisimplicial
category, $\bb S.\bb S.\cc A=D^b(S.S.\cc A)$, and more generally a
multisimplicial category, $\bb S.^n\cc A=D^b(S.^n\cc A)$. There
results a spectrum
   $$n\longmapsto|i.\bb S.^n\cc A|$$
whose structure maps are defined as the map $|i.D^b(\cc
A)|\to\Omega|i.\bb S.\cc A|$ above.

It turns out that the spectrum is a $\Omega$-spectrum beyond the
first term (the additivity theorem is needed to show this, below).
As the spectrum is connective (the $n$th term is
$(n-1)$-connected) an equivalent assertion is that in the sequence
   $$|i.D^b(\cc A)|\lra{}\Omega|i.\bb S.\cc A|\lra{}\Omega\Omega|i.\bb S.\bb S.\cc A|\lra{}\cdots$$
all maps except the first are homotopy equivalences.

Let $\cc A$ be an exact category and let $\cc E$ be its extension
category. There are three natural simplicial maps
$s_*,t_*,q_*:i.\bb S.\cc E\lra{}i.\bb S.\cc A$ induced by
$s,t,q:\cc E\to\cc A$ that take a short exact sequence
   $$\xymatrix{A\ar@{ >->}[r]&C\ar@{->>}[r]&B}$$
to $A$, $C$ and $B$ respectively.

\begin{thm}[Additivity]\label{add}
Let $\cc A$ be an exact category and let $\cc E$ be its extension
category. Then the map
   $$i.\bb S.\cc E\xrightarrow{(s_*,q_*)}i.\bb S.\cc A\times i.\bb S.\cc A$$
is a homotopy equivalence.
\end{thm}

Before proving the theorem we recall the reader certain simplicial
facts.

\begin{lem}[\cite{S}]\label{segal}
Let $X..\to Y..$ be a map of bisimplicial sets. Suppose that for
every $n$, the map $X.{}_n\to Y.{}_n$ is a homotopy equivalence.
Then $X..\to Y..$ is a homotopy equivalence.
\end{lem}

\begin{lem}[\cite{W1}]\label{wald}
Let $X..\to Y..\to Z..$ be a sequence of bisimplicial sets so that
$X..\to Z..$ is constant. Suppose that $X.{}_n\to Y.{}_n\to
Z.{}_n$ is a fibration up to homotopy, for every $n$. We also
require a compatibility with $n$. Suppose further that $Z.{}_n$ is
connected for every $n$. Then $X..\to Y..\to Z..$ is a fibration
up to homotopy.
\end{lem}

Let $\Delta[n]$ denote the simplicial set {\it standard
$n$-simplex},
$\Delta^m\longmapsto\Hom_{\Delta}(\Delta^m,\Delta^n)$. Let
$f:X\lra{}Y$ be a map of simplicial sets and let $y$ be a
$n$-simplex of $Y$. Define a simplicial set $f/(n,y)$ as the
pullback
   $$\begin{CD}
      f/(n,y)@>>>X\\
      @VVV@VVfV\\
      \Delta[n]@>y>>Y.
     \end{CD}$$

\begin{lemb}[\cite{W}]
If for every $u:\Delta^m\lra{}\Delta^n$, and every $y\in Y_n$, the
induced map from $f/(m,u^*y)$ to $f/(n,y)$ is a homotopy
equivalence then for every $(n,y)$ the pullback diagram above is
homotopy cartesian.
\end{lemb}

\begin{lem}\label{szz}
For every $k\geq 0$, the map $f:i_k\bb S.\cc E\to i_k\bb S.\cc A$
sending a string $E$ to the string
$A=A_0\lra{\sim}A_1\lra{\sim}\cdots\lra{\sim}A_k$ satisfies the
hypothesis of Lemma~B.
\end{lem}

By Lemma~B we obtain a homotopy fibration $f/(n,A)\to i_k\bb S.\cc
E\to i_k\bb S.\cc A$ for every simplex $A$ of $i_k\bb S.\cc A$. In
particular the sequence $f/(0,0)\to i_k\bb S.\cc E\to i_k\bb S.\cc
A$ is a homotopy fibration for the unique 0-simplex. The term
$f/(0,0)$ can be identified to the simplicial set $i_k\bb S.\cc
E'$ consisting of the strings $E\in i_k\bb S.\cc E$ such that
$A_j=0$ and $C_j\to B_j$, $j\le k$, is an isomorphism. The latter
simplicial set is homotopy equivalent to $i_k\bb S.\cc A$ by the
sublemma below via the exact equivalence $E\in\cc E'\to B\in\cc
A$. The simplicial set $i_k\bb S.\cc A$ is connected and therefore
the sequence
   $$i.\bb S.\cc A\lra{g}i.\bb S.\cc E\lra{f}i.\bb S.\cc A$$
with $B\bl g\longmapsto 0\rightarrowtail B\twoheadrightarrow B$ is
a fibration by Lemma~\ref{wald}.

Finally, consider a morphism of the latter fibration sequence to
the trivial product fibration sequence,
   $$\begin{CD}
      i.\bb S.\cc A@>>>i.\bb S.\cc E@>>>i.\bb S.\cc A\\
      @V \id VV@VV(s,q)V@VV \id V\\
      i.\bb S.\cc A@>>>i.\bb S.\cc A\times i.\bb S.\cc A@>>>i.\bb S.\cc A.
     \end{CD}$$
The map is a homotopy equivalence on the fibre and on the base,
and hence is so on the total spaces. Thus Lemma~\ref{szz} implies
the additivity theorem.

Let $C$ and $D$ be two simplicial objects in a category $\cc C$
and let $\Delta/\Delta^1$ denote the category of objects over
$\Delta^1$ in $\Delta$; the objects are the maps
$\Delta^n\lra{}\Delta^1$. For any simplicial object $C$ in $\cc C$
let $C^*$ denote the composed functor
   \begin{gather*}
    (\Delta/\Delta^1)^{\op}\lra{}\Delta^{\op}\lra{C}\cc C\\
    (\Delta^n\lra{}\Delta^1)\longmapsto\Delta^n\longmapsto C_n.
   \end{gather*}
Then a {\it simplicial homotopy\/} of maps from $C$ to $D$ is a
natural transformation $C^*\lra{}D^*$~\cite[p.~335]{W}.

\begin{sublem}
Let $\cc A$ and $\cc A'$ be two exact categories. Then an
isomorphism between two exact functors $f,g:\cc A\to\cc A'$
induces a homotopy between $f_*$ and $g_*:i_k\bb S.\cc A\to i_k\bb
S.\cc A'$ for every $k\geq 0$. In particular, every exact
equivalence $\cc A\to\cc A'$ induces a homotopy equivalence
$i_k\bb S.\cc A\to i_k\bb S.\cc A'$.
\end{sublem}

\begin{proof}
The proof is similar to that of~\cite[1.4.1]{W}.
\end{proof}

\renewcommand{\proofname}{Proof of Lemma~\ref{szz}}

\begin{proof}
To simplify the notation the maps $\Ar\Delta^m\lra{}\Ar\Delta^n$
induced by the maps $u:\Delta^m\lra{}\Delta^n$ we denote by the
same letter. We must show that for every $A'\in i_k\bb S_n\cc A$
and $u:\Delta^m\lra{}\Delta^n$ in $\Delta$, the map
$u_*:f/(m,u^*A')\lra{}f/(n,A')$ is a homotopy equivalence. Since
there are maps $v:\Delta^0\lra{}\Delta^n$ and
$w:\Delta^0\lra{}\Delta^m$ such that $uw=v$, it suffices to
consider the special class of maps $\Delta^0\lra{}\Delta^n$.
Indeed, if we proved that both $v_*$ and $w_*$ are homotopy
equivalences, then it would follow that $u_*$ is a homotopy
equivalence, too.

So we must prove the following special case: let $A'$ be a
$n$-simplex of $i_k\bb S.\cc A$, for some $n$, and let
$v_i:\Delta^0\lra{}\Delta^n$ be the map taking $0$ to $i$. Then
for every $i$ the map
   $$v_i{}_*:f/(0,0)\lra{}f/(n,A')$$
is a homotopy equivalence.

A $m$-simplex of $f/(n,A')$ consists of a $m$-simplex $E$ of
$i_k\bb S_m\cc  E$ together with a map
$u:\Ar\Delta^m\lra{}\Ar\Delta^n$ such that $u^*A'=E_{(0,0)}$. The
map sending $E$ to $E_{(1,1)}$ induces a map
$p:f/(n,A')\lra{}i_k\bb S.\cc A$. It will suffice to show that $p$
is a homotopy equivalence. Indeed, $p$ is left inverse to the
composed map
   $$i_k\bb S.\cc A\lra{\beta_*}f/(0,0)\lra{v_i{}_*}f/(n,A),$$
therefore if $p$ is a homotopy equivalence then so is
$v_i{}_*\beta$ and hence also $v_i{}_*$, since the map $\beta$
taking $B$ to $0\rightarrowtail B\twoheadrightarrow B$ is a
homotopy equivalence by the sublemma above. This implies $v_i{}_*$
is a homotopy equivalence, too. To prove that $p$ is a homotopy
equivalence, it suffices to show that the particular map $v_n{}_*$
is a homotopy equivalence, because $pv_n{}_*\beta=1$.

We shall construct the homotopy by lifting the simplicial homotopy
that contracts $\Delta[n]$ to its last vertex. This simplicial
homotopy is given by a map of the composed functors
   \begin{gather*}
    (\Delta/\Delta^1)^{\op}\lra{}\Delta^{\op}\lra{}\Sets\\
    (\Delta^m\lra{}\Delta^1)\longmapsto\Delta^m\longmapsto\Hom_\Delta(\Delta^m,\Delta^n)
   \end{gather*}
to itself. Precisely, the functor takes $v:\Delta^m\lra{}\Delta^1$
to
   $$(u:\Delta^m\lra{}\Delta^n)\longmapsto(\bar u:\Delta^m\lra{}\Delta^n)$$
where $\bar u$ is defined as the composite
   $$\Delta^m\xrightarrow{(u,v)}\Delta^n\times\Delta^1\lra{w}\Delta^n$$
and where $w(j,0)=j$ and $w(j,1)=n$.

A lifting of this homotopy to one on $f/(n,A')$ is a map taking
$v:\Delta^m\lra{}\Delta^1$ to
   $$(E,u)\to(\bar E,\bar u)$$
with $\bar E_{(0,0)}=\bar u^*A'$. We shall depict the elements of
$i_k\bb S_m\cc E$ as diagrams
   $$\xymatrix{
     A_0\ar @{ >->} [d] \ar[r]^{a_1}&A_1\ar@{ >->}[d]\ar[r]^{a_2}&\cdots\ar[r]^{a_k}&A_k\ar@{ >->}[d]\\
     C_0\ar@{->>}[d]\ar[r]^{c_1}&C_1\ar@{->>}[d]\ar[r]^{c_2}&\cdots\ar[r]^{c_k}&C_k\ar@{->>}[d]\\
     B_0\ar[r]^{b_1}&B_1\ar[r]^{b_2}&\cdots\ar[r]^{b_k}&B_k
     }$$
with $A_i=u^*A'_i$, $a_i=u^*(a_i':A_{i-1}'\to A_i')$, and
$(a_i,c_i,b_i):E_{i-1}\to E_i$ isomorphisms in $D^b(S_m\cc E)$,
$i\le k$. Each vertical map is represented by a commutative
diagram
  $$\xymatrix{ A_{i-1} \ar@{ >->}[r]& C_{i-1} \ar@{->>}[r] & B_{i-1}\\
  Y_i \ar@{ >->}[r]\ar[u]^{s_i}\ar[d]_{t_i} & Z_i\ar@{->>}[r]\ar[u]_{s'_i}\ar[d]^{t'_i} & W_i\ar[u]_{s''_i}\ar[d]^{t''_i}\\
  A_i \ar@{ >->}[r] & C_i\ar@{->>}[r] & B_i}$$
in $C^b(S_m\cc A)$ with the vertical maps quasi-isomorphisms and
$t_is_i^{-1}=a_i$.

Since $u\le\bar u$ by construction, it follows that there is a
bimorphism $\phi:u\to\bar u$. This bimorphism is actually unique,
because we deal with maps of posets. This yields a map
$\phi^*_{A_i}:A_i\to\bar A_i$ for every $i$ where $\bar A_i=\bar
u^*A_i'$. By assumption, each morphism $a_i:A_{i-1}\to A_i$ equals
to $u^*(a_i')$ where $a_i':A_{i-1}'\to A_i'$ is an isomorphism in
$D^b(S_n\cc A)$ represented by the equivalence class of a diagram
   $$A_{i-1}'\bl{p_i'}\longleftarrow X_i'\lra{q_i'}A_i'$$
with $p_i', q_i'$ quasi-isomorphisms in $C^b(S_n\cc A)$. Then
$a_i=q_ip_i^{-1}$ where $p_i=u^*(p_i')$ and $q_i=u^*(q_i')$.

There is a common denominator
   \begin{equation}\label{den}
    \xymatrix{&&V_i\ar[dl]_{u_i}\ar[dr]^{v_i}\\
      &Y_i\ar[dl]_{s_i}\ar[drrr]^(.2){t_i}&& X_i\ar[dlll]_(.2){p_i}\ar[dr]^{q_i}\\
      A_{i-1}&&&&A_i}
   \end{equation}
By Lemma~\ref{imp} there exists a complex $U_i\in C^b(S_m\cc A)$
and a commutative diagram
   $$\xymatrix{
      C_{i-1} & U_i\ar[l]_{\sigma_i}\ar[r]^{\tau_i} & C_i\\
      A_{i-1}\ar@{ >->}[u] & V_i \ar[l]^{p_iv_i}\ar@{ >->}[u]^{x_i}\ar[r]_{q_iv_i}&A_i\ar@{ >->}[u]}$$
representing the same morphism in $D^b([S_m\cc A]^{\Delta^1})$
   \begin{equation}\label{trang}
    \xymatrix{
      C_{i-1} & Z_i\ar[l]_{s'_i}\ar[r]^{t'_i} & C_i\\
      A_{i-1}\ar@{ >->}[u] & Y_i \ar[l]^{s_i}\ar@{ >->}[u]\ar[r]_{t_i}&A_i\ar@{ >->}[u]}
   \end{equation}

Let $\bar a_i:\bar A_{i-1}\to\bar A_i$ be the map represented by
the equivalence class of the diagram
   $$\bar A_{i-1}\bl{\bar p_i}\longleftarrow\bar X_i=\bar u^*X_i'\lra{\bar q_i}\bar A_i$$
with $\bar p_i=\bar u^*(p_i'),\bar q_i=\bar u^*(q_i')$. Then $\bar
a_i$ is an isomorphism since $a_i'$ is so.

We obtain a commutative diagram
   \begin{equation}\label{33}
    \xymatrix{
      C_{i-1} & U_i\ar[l]_{\sigma_i}\ar[r]^{\tau_i} & C_i\\
      A_{i-1}\ar@{ >->}[u]\ar[d]_{\phi^*_{A_{i-1}}}
      &V_i\ar[l]^{p_iv_i}\ar@{ >->}[u]^{x_i}\ar[d]^{\phi^*_{X_i}v_i}\ar[r]_{q_iv_i}&A_i\ar@{ >->}[u]\ar[d]^{\phi^*_{A_i}}\\
      \bar A_{i-1}&\bar X_i\ar[l]^{\bar p_i}\ar[r]_{\bar q_i}&\bar A_i}
   \end{equation}
giving an isomorphism $(\bar a_i,a_i,c_i)$ in $D^b(\wt{[S_m\cc
A]\ulcorner})$.

\begin{sublem}

The map $(\bar a_i,a_i,c_i)$ represented by diagram~\eqref{33} is
well defined that is it does not depend on:

$(1)$ the choice of a common denominator~\eqref{den};

$(2)$ the choice of a representative for $(a_i,c_i,b_i):E_{i-1}\to
E_i$;

$(3)$ the choice of a representative for $a'_i:A'_{i-1}\to A'_i$.

\end{sublem}

\renewcommand{\proofname}{Proof}

\begin{proof}
Let us check~(1). Suppose we are given the following diagram.
    $$\xymatrix{&&V_i\ar[dl]_{u_i}\ar[drr]^{v_i}&V_i'\ar[dll]_{w_i}\ar[dr]^{z_i}\\
      &Y_i\ar[dl]_{s_i}\ar[drrrr]^(.2){t_i}&&& X_i\ar[dllll]_(.2){p_i}\ar[dr]^{q_i}\\
      A_{i-1}&&&&&A_i}$$
We have to show that
   \begin{equation}\label{int}
    \xymatrix{
      C_{i-1} & U_i\ar[l]_{\sigma_i}\ar[r]^{\tau_i} & C_i\\
      A_{i-1}\ar@{ >->}[u]\ar[d]_{\phi^*_{A_{i-1}}}
      &V_i\ar[l]^{p_iv_i}\ar@{ >->}[u]^{x_i}\ar[d]^{\phi^*_{X_i}v_i}\ar[r]_{q_iv_i}&A_i\ar@{ >->}[u]\ar[d]^{\phi^*_{A_i}}\\
      \bar A_{i-1}&\bar X_i\ar[l]^{\bar p_i}\ar[r]_{\bar
      q_i}&\bar A_i}
      \textrm{is equivalent to}
      \xymatrix{
      C_{i-1} & U'_i\ar[l]_{\sigma'_i}\ar[r]^{\tau'_i} & C_i\\
      A_{i-1}\ar@{ >->}[u]\ar[d]_{\phi^*_{A_{i-1}}}
      &V'_i\ar[l]^{p_iz_i}\ar@{ >->}[u]^{x'_i}\ar[d]^{\phi^*_{X_i}z_i}\ar[r]_{q_iz_i}&A_i\ar@{ >->}[u]\ar[d]^{\phi^*_{A_i}}\\
      \bar A_{i-1}&\bar X_i\ar[l]^{\bar p_i}\ar[r]_{\bar q_i}&\bar A_i}
   \end{equation}
By Lemma~\ref{imp} there is a common denominator
   $$
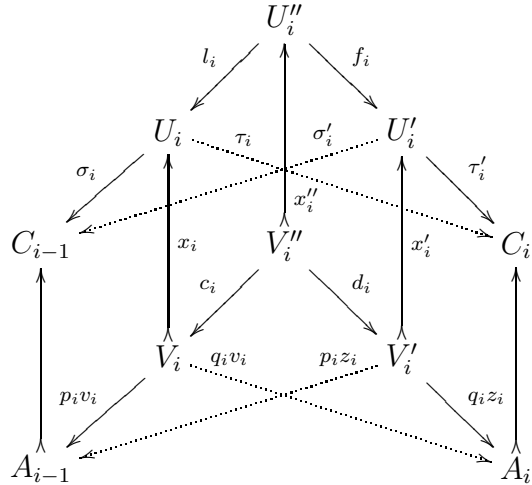
\xymatrix{
      &&U_i''\ar[dl]_{l_i}\ar[dr]^{f_i}\\
      &U_i\ar[dl]_{\sigma_i}\ar@{.>}[drrr]^(.2){\tau_i}&& U'_i\ar@{.>}[dlll]_(.2){\sigma'_i}\ar[dr]^{\tau'_i}\\
      C_{i-1}&&V_i''\ar[dl]_{c_i}\ar@{ >->}[uu]_(.2){x_i''}\ar[dr]^{d_i}&&C_i\\
      &V_i\ar[dl]_{p_iv_i}\ar@{ >->}[uu]_{x_i}\ar@{.>}[drrr]^(.15){q_iv_i}&&
      V'_i\ar@{.>}[dlll]_(.15){p_iz_i}\ar@{ >->}[uu]_{x_i'}\ar[dr]^{q_iz_i}\\
      A_{i-1}\ar@{ >->}[uu]&&&&A_i\ar@{ >->}[uu]}$$
Fix homotopies
$(e^n,g^n):(V''^{n}_i\lra{x_i''}U''^{n}_i)\to(A^{n-1}_{i-1}\lra{}C^{n-1}_{i-1})$
and
$(h^n,m^n):(V''^{n}_i\lra{x_i''}U''^{n}_i)\to(A^{n-1}_{i}\lra{}C^{n-1}_{i})$
for $(p_iv_ic_i,\sigma_il_i)\sim(p_iz_id_i,\sigma_i'f_i)$ and
$(q_iv_ic_i,\tau_il_i)\sim(q_iz_id_i,\tau_i'f_i)$ respectively.
The latter diagram fits into the diagram
   
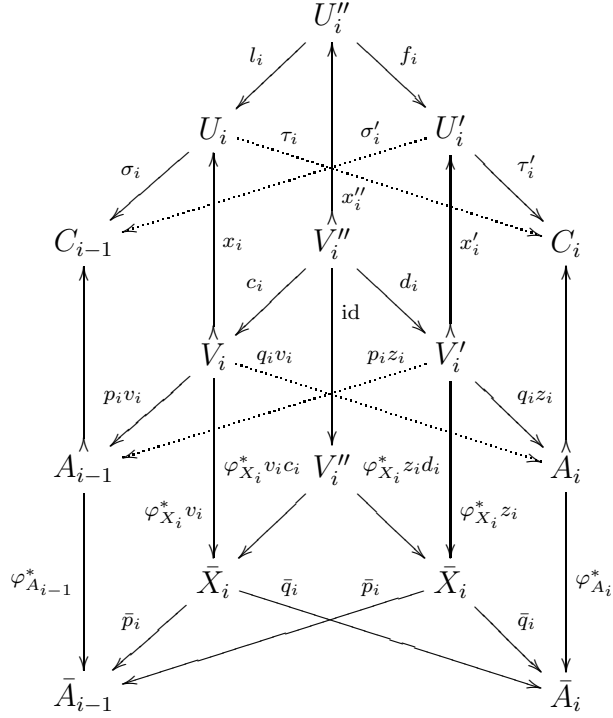
\begin{equation}\label{fra}
    \xymatrix{
      &&U_i''\ar[dl]_{l_i}\ar[dr]^{f_i}\\
      &U_i\ar[dl]_{\sigma_i}\ar@{.>}[drrr]^(.2){\tau_i}&& U'_i\ar@{.>}[dlll]_(.2){\sigma'_i}\ar[dr]^{\tau'_i}\\
      C_{i-1}&&V_i''\ar[dl]_{c_i}\ar@{ >->}[uu]_(.2){x_i''}\ar[dd]^(.3){\id}\ar[dr]^{d_i}&&C_i\\
      &V_i\ar[dl]_{p_iv_i}\ar@{ >->}[uu]_{x_i}\ar[dd]_(.7){\phi^*_{X_i}v_i}\ar@{.>}[drrr]^(.15){q_iv_i}&&
      V'_i\ar@{.>}[dlll]_(.15){p_iz_i}\ar@{ >->}[uu]_{x_i'}\ar[dd]^(.7){\phi^*_{X_i}z_i}\ar[dr]^{q_iz_i}\\
      A_{i-1}\ar@{ >->}[uu]\ar[dd]_{\phi^*_{A_{i-1}}}&&
      V_i''\ar[dl]_(.2){\phi^*_{X_i}v_ic_i}\ar[dr]^(.2){\phi^*_{X_i}z_id_i}&&A_i\ar@{ >->}[uu]\ar[dd]^{\phi^*_{A_i}}\\
      &\bar X_i\ar[dl]_{\bar p_i}\ar[drrr]^(.2){\bar q_i}&& \bar X_i\ar[dlll]_(.2){\bar p_i}\ar[dr]^{\bar q_i}\\
      \bar A_{i-1}&&&&\bar A_i}
   \end{equation}

We want to show that this diagram is a common denominator
for~\eqref{int}. Put $k^n=\phi^{*n-1}_{A_{i-1}}\circ
e^n:V''^n\to\bar A^{n-1}_{i-1}$. Then
   \begin{gather*}
    \bar p_i\phi^*_{X_i}v_ic_i-\bar
    p_i\phi^*_{X_i}z_id_i=\phi^*_{A_{i-1}}p_iv_ic_i-\phi^*_{A_{i-1}}p_iz_id_i=\\
     =\phi^*_{A_{i-1}}e\partial+\phi^*_{A_{i-1}}\partial e=k\partial+\partial\phi^*_{A_{i-1}}e=k\partial+\partial k.
   \end{gather*}
This shows that $(\sigma_il_i,p_iv_ic_i,\bar
p_i\phi^*_{X_i}v_ic_i)\bl{(g,e,k)}\sim(\sigma'_if_i,p_iz_id_i,\bar
p_i\phi^*_{X_i}z_id_i)$. A homotopy between
$(\tau_il_i,q_iv_ic_i,\bar q_i\phi^*_{X_i}v_ic_i)$ and
$(\tau'_if_i,q_iz_id_i,\bar q_i\phi^*_{X_i}z_id_i)$ is similarly
checked.

Let us check~(2). Suppose
   $$\xymatrix{C_{i-1} & Z'_i\ar[l]_{\alpha'_i}\ar[r]^{\beta'_i} & C_i\\
      A_{i-1}\ar@{ >->}[u] & Y'_i
      \ar[l]^{\alpha_i}\ar@{ >->}[u]\ar[r]_{\beta_i}&A_i\ar@{ >->}[u]}$$
is equivalent to~\eqref{trang}. There is a homotopy commutative
diagram
   $$\xymatrix{
      &V_i\ar[dl]_{u_i}\ar[dr]^{v_i}&&V_i'\ar[dl]_{z_i}\ar[dr]^{w_i}\\
      Y_i\ar[d]_{s_i}\ar[drrrr]^(.2){t_i}&&X_i\ar[dll]_(0.2){p_i}\ar[drr]^(0.2){q_i}
      &&Y_i'\ar[dllll]_(0.2){\alpha_i}\ar[d]^{\beta_i}\\
      A_{i-1}&&&&A_i
     }$$

We have to show that
   \begin{equation}\label{hu}
    \xymatrix{
      C_{i-1} & U_i\ar[l]_{\sigma_i}\ar[r]^{\tau_i} & C_i\\
      A_{i-1}\ar@{ >->}[u]\ar[d]_{\phi^*_{A_{i-1}}}
      &V_i\ar[l]^{p_iv_i}\ar@{ >->}[u]^{x_i}\ar[d]^{\phi^*_{X_i}v_i}\ar[r]_{q_iv_i}&A_i\ar@{ >->}[u]\ar[d]^{\phi^*_{A_i}}\\
      \bar A_{i-1}&\bar X_i\ar[l]^{\bar p_i}\ar[r]_{\bar
      q_i}&\bar A_i}
      \textrm{is equivalent to}
      \xymatrix{
      C_{i-1} & U'_i\ar[l]_{\sigma'_i}\ar[r]^{\tau'_i} & C_i\\
      A_{i-1}\ar@{ >->}[u]\ar[d]_{\phi^*_{A_{i-1}}}
      &V'_i\ar[l]^{p_iz_i}\ar@{ >->}[u]^{x'_i}\ar[d]^{\phi^*_{X_i}z_i}\ar[r]_{q_iz_i}&A_i\ar@{ >->}[u]\ar[d]^{\phi^*_{A_i}}\\
      \bar A_{i-1}&\bar X_i\ar[l]^{\bar p_i}\ar[r]_{\bar q_i}&\bar A_i}
   \end{equation}
It follows from Lemma~\ref{imp} and our assumption that
   $$\xymatrix{C_{i-1} & U_i\ar[l]_{\sigma_i}\ar[r]^{\tau_i} & C_i\\
      A_{i-1}\ar@{ >->}[u] & V_i
      \ar[l]^{p_iv_i}\ar@{ >->}[u]^{x_i}\ar[r]_{q_iv_i}&A_i\ar@{ >->}[u]}
      \textrm{ is equivalent to }
      \xymatrix{C_{i-1} & Z_i\ar[l]\ar[r] & C_i\\
      A_{i-1}\ar@{ >->}[u] & Y_i
      \ar[l]^{s_i}\ar@{ >->}[u]\ar[r]_{t_i}&A_i\ar@{ >->}[u]}
      \textrm{ is equivalent to }$$

   $$\xymatrix{C_{i-1} & Z'_i\ar[l]\ar[r] & C_i\\
      A_{i-1}\ar@{ >->}[u] & Y'_i
      \ar[l]^{\alpha_i}\ar@{ >->}[u]\ar[r]_{\beta_i}&A_i\ar@{ >->}[u]}
      \textrm{ is equivalent to }
      \xymatrix{C_{i-1} & U'_i\ar[l]\ar[r] & C_i\\
      A_{i-1}\ar@{ >->}[u] & V'_i
      \ar[l]^{p_iz_i}\ar@{ >->}[u]^{x'_i}\ar[r]_{q_iz_i}&A_i\ar@{ >->}[u]}$$
One can now construct the diagram~\eqref{fra} yielding a common
denominator for~\eqref{hu}. This implies~(2). It remains to
check~(3).

Let $A_{i-1}'\bl{r_i'}\longleftarrow W_i'\lra{n_i'}A_i'$ be
another representative for $a_i'$. There is a homotopy commutative
diagram
   $$\xymatrix{
      &V_i\ar[dl]_{v_i}\ar[dr]^{u_i}&&V_i'\ar[dl]_{w_i}\ar[dr]^{z_i}\\
      X_i\ar[d]_{p_i}\ar[drrrr]^(.2){q_i}&&Y_i\ar[dll]_(0.2){s_i}\ar[drr]^(0.2){t_i}
      &&W_i\ar[dllll]_(0.2){r_i}\ar[d]^{n_i}\\
      A_{i-1}&&&&A_i
     }$$
with $r_i=u^*(r_i'),n_i=u^*(n_i')$. We have to show that
   $$\xymatrix{
      C_{i-1} & U_i\ar[l]_{\sigma_i}\ar[r]^{\tau_i} & C_i\\
      A_{i-1}\ar@{ >->}[u]_{\hskip0.6cm\textrm{(I)}}\ar[d]_{\phi^*_{A_{i-1}}}
      &V_i\ar[l]^{p_iv_i}\ar@{ >->}[u]^{x_i}\ar[d]^{\phi^*_{X_i}v_i}\ar[r]_{q_iv_i}&A_i\ar@{ >->}[u]\ar[d]^{\phi^*_{A_i}}\\
      \bar A_{i-1}&\bar X_i\ar[l]^{\bar p_i}\ar[r]_{\bar
      q_i}&\bar A_i}
      \textrm{is equivalent to}
      \xymatrix{
      C_{i-1} & U'_i\ar[l]_{\sigma'_i}\ar[r]^{\tau'_i} & C_i\\
      A_{i-1}\ar@{ >->}[u]_{\hskip0.5cm\textrm{(II)}}\ar[d]_{\phi^*_{A_{i-1}}}
      &V'_i\ar[l]^{r_iz_i}\ar@{ >->}[u]^{x'_i}\ar[d]^{\phi^*_{X_i}z_i}\ar[r]_{n_iz_i}&A_i\ar@{ >->}[u]\ar[d]^{\phi^*_{A_i}}\\
      \bar A_{i-1}&\bar X_i\ar[l]^{\bar r_i}\ar[r]_{\bar n_i}&\bar A_i}$$

There are common denominators in $C^b(S_n\cc A)$ and in
$C^b(S_m\cc A)$, respectively.
   $$
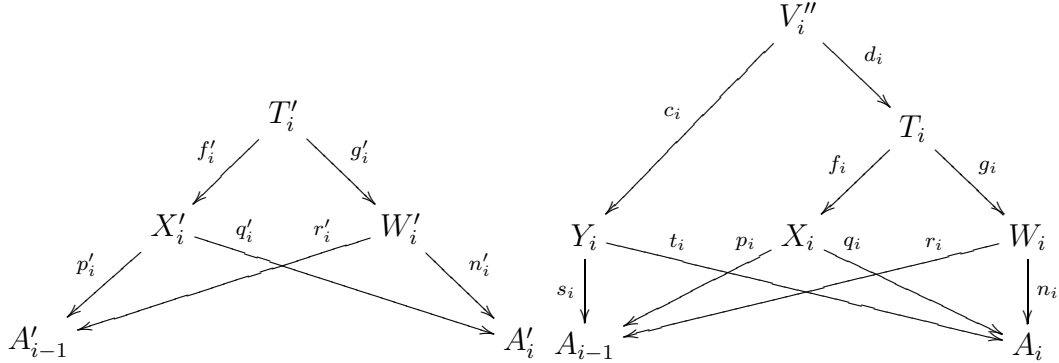
\xymatrix{{}\\ &&T'_i\ar[dl]_{f'_i}\ar[dr]^{g'_i}\\
      &X'_i\ar[dl]_{p'_i}\ar[drrr]^(.2){q'_i}&& W'_i\ar[dlll]_(.2){r'_i}\ar[dr]^{n'_i}\\
      A_{i-1}'&&&&A_i'}
     \xymatrix{
      &&V_i''\ar[ddll]_{c_i}\ar[dr]^{d_i}\\ &&&T_i\ar[dl]_{f_i}\ar[dr]^{g_i}\\
      Y_i\ar[d]_{s_i}\ar[drrrr]^(.2){t_i}&&X_i\ar[dll]_(0.2){p_i}\ar[drr]^(0.2){q_i}
      &&W_i\ar[dllll]_(0.2){r_i}\ar[d]^{n_i}\\
      A_{i-1}&&&&A_i
     }$$
It follows from~(1) that
   $$\textrm{(I)} \sim \xymatrix{
      C_{i-1} & E_i\ar[l]_{\pi_i}\ar[r]^{\rho_i} & C_i\\
      A_{i-1}\ar@{ >->}[u]_{\hskip0.5cm\textrm{(III)}}\ar[d]_{\phi^*_{A_{i-1}}}
      &V_i''\ar[l]^{p_if_id_i}\ar@{ >->}[u]^{\kappa_i}\ar[d]^{\phi^*_{X_i}f_id_i}\ar[r]_{q_if_id_i}&A_i\ar@{ >->}[u]\ar[d]^{\phi^*_{A_i}}\\
      \bar A_{i-1}&\bar X_i\ar[l]^{\bar p_i}\ar[r]_{\bar
      q_i}&\bar A_i}
      \textrm{and (II) $\sim$ }
      \xymatrix{
      C_{i-1} & E'_i\ar[l]_{\pi'_i}\ar[r]^{\rho'_i} & C_i\\
      A_{i-1}\ar@{ >->}[u]_{\hskip0.5cm\textrm{(IV)}}\ar[d]_{\phi^*_{A_{i-1}}}
      &V_i''\ar[l]^{r_ig_id_i}\ar@{ >->}[u]^{\kappa'_i}\ar[d]^{\phi^*_{W_i}g_id_i}\ar[r]_{n_ig_id_i}&A_i\ar@{ >->}[u]\ar[d]^{\phi^*_{A_i}}\\
      \bar A_{i-1}&\bar W_i\ar[l]^{\bar r_i}\ar[r]_{\bar n_i}&\bar A_i}$$
Since $\phi^*_{X_i}f_i=\bar f_i\phi^*_{T_i}$ and
$\phi^*_{W_i}g_i=\bar g_i\phi^*_{T_i}$, it follows that
   $$\textrm{(III)} \sim \xymatrix{
      C_{i-1} & E_i\ar[l]_{\pi_i}\ar[r]^{\rho_i} & C_i\\
      A_{i-1}\ar@{ >->}[u]_{\hskip0.6cm\textrm{(V)}}\ar[d]_{\phi^*_{A_{i-1}}}
      &V_i''\ar[l]^{p_if_id_i}\ar@{ >->}[u]^{\kappa_i}\ar[d]^{\phi^*_{T_i}d_i}\ar[r]_{q_if_id_i}&A_i\ar@{ >->}[u]\ar[d]^{\phi^*_{A_i}}\\
      \bar A_{i-1}&\bar T_i\ar[l]^{\bar p_i\bar f_i}\ar[r]_{\bar q_i\bar f_i}&\bar A_i}
      \textrm{and (IV) $\sim$ }
      \xymatrix{
      C_{i-1} & E'_i\ar[l]_{\pi'_i}\ar[r]^{\rho'_i} & C_i\\
      A_{i-1}\ar@{ >->}[u]_{\hskip0.5cm\textrm{(VI)}}\ar[d]_{\phi^*_{A_{i-1}}}
      &V_i''\ar[l]^{r_ig_id_i}\ar@{ >->}[u]^{\kappa'_i}\ar[d]^{\phi^*_{T_i}d_i}\ar[r]_{n_ig_id_i}&A_i\ar@{ >->}[u]\ar[d]^{\phi^*_{A_i}}\\
      \bar A_{i-1}&\bar T_i\ar[l]^{\bar r_i\bar g_i}\ar[r]_{\bar n_i\bar g_i}&\bar A_i}$$
Fix homotopies $p_i'f_i'\bl{k'}\sim r_i'g_i', n_i'g_i'\bl{e'}\sim
q_i'f_i'$. Then $(k,e)=(u^*(k'),u^*(e'))$ and $(\bar k,\bar
e)=(\bar u^*(k'),\bar u^*(e'))$ give homotopies $p_if_i\bl{k}\sim
r_ig_i, n_ig_i\bl{e}\sim q_if_i$ and $\bar p_i\bar f_i\bl{\bar
k}\sim \bar r_i\bar g_i, \bar n_i\bar g_i\bl{\bar e}\sim \bar
q_i\bar f_i$. The proof of Lemma~\ref{imp} shows that~(V) can be
embedded into the diagram
   $$
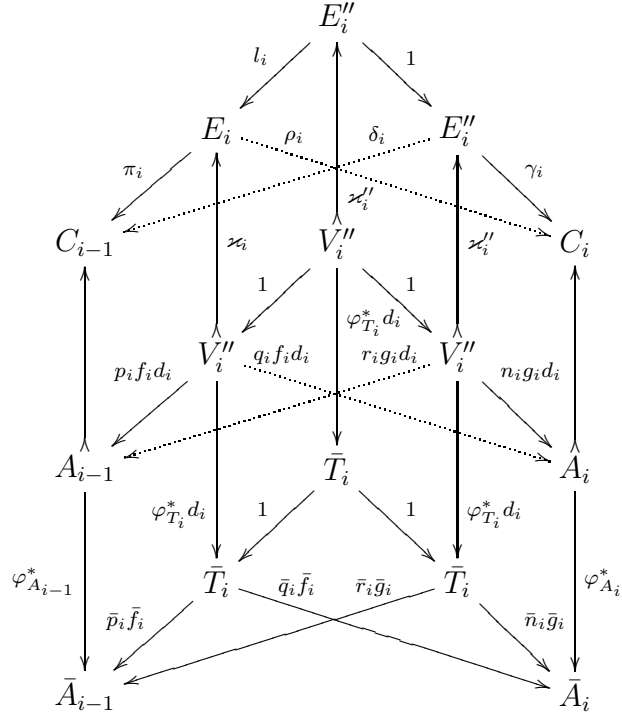
\xymatrix{
      &&E_i''\ar[dl]_{l_i}\ar[dr]^{1}\\
      &E_i\ar[dl]_{\pi_i}\ar@{.>}[drrr]^(.2){\rho_i}&& E''_i\ar@{.>}[dlll]_(.2){\delta_i}\ar[dr]^{\gamma_i}\\
      C_{i-1}&&V_i''\ar[dl]_{1}\ar@{ >->}[uu]_(.2){\kappa''_i}\ar[dd]^(.35){\phi^*_{T_i}d_i}\ar[dr]^{1}&&C_i\\
      &V''_i\ar[dl]_(.3){p_if_id_i}\ar@{ >->}[uu]_{\kappa_i}\ar[dd]_(.7){\phi^*_{T_i}d_i}\ar@{.>}[drrr]^(.15){q_if_id_i}&&
      V''_i\ar@{.>}[dlll]_(.15){r_ig_id_i}\ar@{ >->}[uu]_{\kappa''_i}\ar[dd]^(.7){\phi^*_{T_i}d_i}\ar[dr]^(.3){n_ig_id_i}\\
      A_{i-1}\ar@{ >->}[uu]\ar[dd]_{\phi^*_{A_{i-1}}}&&
      \bar T_i\ar[dl]_{1}\ar[dr]^1&&A_i\ar@{ >->}[uu]\ar[dd]^{\phi^*_{A_i}}\\
      &\bar T_i\ar[dl]_{\bar p_i\bar f_i}\ar[drrr]^(.2){\bar q_i\bar f_i}
      && \bar T_i\ar[dlll]_(.2){\bar r_i\bar g_i}\ar[dr]^{\bar n_i\bar g_i}\\
      \bar A_{i-1}&&&&\bar A_i}$$
in such a way that the upper ``roof" is a common denominator with
homotopies
$(p_if_id_i,\pi_il_i)\bl{(kd_i,h)}\sim(r_ig_id_i,\delta_i)$ and
$(n_ig_id_i,\gamma_i)\bl{(ed_i,\ell)}\sim(q_if_id_i,\rho_il_i)$.
We have homotopies $(\bar p_i\bar f_i,p_if_id_i,\pi_il_i)\bl{(\bar
k,kd_i,h)}\sim(\bar r_i\bar g_i,r_ig_id_i,\delta_i)$ and $(\bar
n_i\bar g_i,n_ig_id_i,\gamma_i)\bl{(\bar e,ed_i,\ell)}\sim(\bar
q_i\bar f_i,q_if_id_i,\rho_il_i)$ showing that
   $$\xymatrix{
      C_{i-1} & E_i\ar[l]_{\pi_i}\ar[r]^{\rho_i} & C_i\\
      A_{i-1}\ar@{ >->}[u]_{\hskip0.6cm\textrm{(V)}}\ar[d]_{\phi^*_{A_{i-1}}}
      &V_i''\ar[l]^{p_if_id_i}\ar@{ >->}[u]^{\kappa_i}\ar[d]^{\phi^*_{T_i}d_i}\ar[r]_{q_if_id_i}&A_i\ar@{ >->}[u]\ar[d]^{\phi^*_{A_i}}\\
      \bar A_{i-1}&\bar T_i\ar[l]^{\bar p_i\bar f_i}\ar[r]_{\bar q_i\bar f_i}&\bar A_i}
      \textrm{ is equivalent to }
      \xymatrix{
      C_{i-1} & E''_i\ar[l]_{\delta_i}\ar[r]^{\gamma_i} & C_i\\
      A_{i-1}\ar@{ >->}[u]_{\hskip0.4cm\textrm{(VII)}}\ar[d]_{\phi^*_{A_{i-1}}}
      &V_i''\ar[l]^{r_ig_id_i}\ar@{ >->}[u]^{\kappa''_i}\ar[d]^{\phi^*_{T_i}d_i}\ar[r]_{n_ig_id_i}&A_i\ar@{ >->}[u]\ar[d]^{\phi^*_{A_i}}\\
      \bar A_{i-1}&\bar T_i\ar[l]^{\bar r_i\bar g_i}\ar[r]_{\bar n_i\bar g_i}&\bar A_i}$$

The proof of (1) shows that $\textrm{(VI)}\sim\textrm{(VII)}$,
hence
$\textrm{(I)}\sim\textrm{(III)}\sim\textrm{(V)}\sim\textrm{(VII)}\sim\textrm{(VI)}\sim\textrm{(IV)}
\sim\textrm{(II)}$ as required.
\end{proof}

For every $i\le k$, we can construct a pushout diagram in
$C^b(S_m\cc A)$
   \begin{equation}\label{push}
    \xymatrix{
     A_i\ar@{ >->}[r]\ar[d]_{\phi^*} & C_i\ar@{->>}[r]\ar[d] & B_i\ar@{=}[d]\\
     \bar A_i\ar@{ >->}[r] &\bar C_i\ar@{->>}[r] &\bar B_i}
   \end{equation}
The diagrams~\eqref{33} give a string of isomorphisms in
$D^b(\wt{[S_m\cc A]\ulcorner})$
   $$\xymatrix{
     C_0\ar[r]^{c_1} &C_1\ar[r]^{c_2}&\cdots\ar[r]^{c_k}&C_k\\
     A_0\ar@{ >->}[u]\ar[d]_{\phi^*}\ar[r]^{a_1}&A_1\ar@{ >->}[u]\ar[d]_{\phi^*}\ar[r]^{a_2}&\cdots\ar[r]^{a_k}
     &A_k\ar@{ >->}[u]\ar[d]^{\phi^*}\\
     \bar A_0\ar[r]^{\bar a_1}&\bar A_1\ar[r]^{\bar a_2}&\cdots\ar[r]^{\bar a_k}&\bar A_k\\
     }$$
By Corollary~\ref{van} one can fit it into a string of
isomorphisms in $D^b(\wt{[S_m\cc A]^\Box})$
   $$\xymatrix@!0{
     &C_0\ar[rr]^{c_1}\ar'[d][dd] && C_1\ar'[d][dd] \ar[rr]^{c_2}\ar'[d][dd] && C_2\ar'[d][dd] \ar[rr]^{c_3}\ar'[d][dd]
     && \cdots \ar[rr]^{c_k} && C_k\ar[dd]\\
     A_0\ar@{>->}[ur]\ar[rr]\ar[dd] && A_1\ar@{>->}[ur]\ar[dd] \ar[ur]\ar[rr]\ar[dd]
     && A_2\ar@{>->}[ur]\ar[dd] \ar[ur]\ar[rr]\ar[dd]&& \cdots\ar[rr] &&A_k\ar@{>->}[ur]\ar[dd]\\
     &\bar C_0\ar@{.>}'[r]_{\bar c_1}[rr]&&\bar C_1\ar@{.>}'[r]_{\bar c_2}[rr] &&\bar C_2\ar@{.>}_{\bar c_3}[rr]
     && \cdots \ar@{.>}'[r]_{\bar c_k}[rr] &&\bar C_k\\
     \bar A_0\ar[rr]_{\bar a_1}\ar@{>->}[ur]&&\bar A_1\ar@{>->}[ur]\ar[rr]_{\bar a_2}&&\bar A_2\ar@{>->}[ur]\ar[rr]_{\bar a_3}
     && \cdots \ar[rr]_{\bar a_k} && \bar A_k\ar@{>->}[ur]
     }$$

Finally, we consider the diagram constructed as above
   $$\xymatrix@!0{
     &\bar C_0\ar[rr]^{\bar c_1}\ar@{->>}'[d][dd] &&\bar C_1\ar@{->>}'[d][dd] \ar[rr]^{\bar c_2}\ar@{->>}'[d][dd]
     &&\bar C_2\ar@{->>}'[d][dd] \ar[rr]^{\bar c_3}\ar@{->>}'[d][dd] && \cdots\ar[rr]^{\bar c_k} &&\bar C_k\ar@{->>}[dd]\\
     \bar A_0\ar@{>->}[ur]\ar[rr]\ar[dd] &&\bar A_1\ar@{>->}[ur]\ar[dd] \ar[ur]\ar[rr]\ar[dd]
     &&\bar A_2\ar@{>->}[ur]\ar[dd] \ar[ur]\ar[rr]\ar[dd]&& \cdots\ar[rr] &&\bar A_k\ar@{>->}[ur]\ar[dd]\\
     &\bar B_0\ar@{.>}'[r]_{\bar b_1}[rr]&&\bar B_1\ar@{.>}'[r]_{\bar b_2}[rr] &&\bar B_2\ar@{.>}_{\bar b_3}[rr]
     && \cdots \ar@{.>}'[r]_{\bar b_k}[rr] &&\bar B_k\\
     0\ar[rr]\ar[ur]&&0\ar[ur]\ar[rr]&&0\ar[ur]\ar[rr]&& \cdots\ar[rr]&& 0\ar[ur]
     }$$

Construction of the simplex $(\bar E,\bar u)$ is completed. It is
given by the diagram
   $$\xymatrix{
     \bar A_0\ar@{ >->}[d]\ar[r]^{\bar a_1}&\bar A_1\ar@{ >->}[d]\ar[r]^{\bar a_2}&\cdots\ar[r]^{\bar a_k}&\bar A_k\ar@{ >->}[d]\\
     \bar C_0\ar@{->>}[d]\ar[r]^{\bar c_1}&\bar C_1\ar@{->>}[d]\ar[r]^{\bar c_2}&\cdots\ar[r]^{\bar c_k}&\bar C_k\ar@{->>}[d]\\
     \bar B_0\ar[r]^{\bar b_1}&\bar B_1\ar[r]^{\bar b_2}&\cdots\ar[r]^{\bar b_k}&\bar B_k
     }$$
with $\bar A_i=\bar u^*A'_i$ and $(\bar a_i,\bar c_i,\bar
b_i):\bar E_{i-1}\to\bar E_i$ isomorphisms in $D^b(S_m\cc E)$,
$i\le k$.

We have to verify that the construction of the simplex $(\bar
E,\bar u)$ is compatible with the structure maps of the category
$\Delta/\Delta^1$; that is, if we replace $\Delta^m$ by
$\Delta^{m'}$ throughout, by means of some map
$\Delta^{m'}\to\Delta^m$, then the strucure map in $i_k\bb S.\cc
E$ induced by $\Delta^{m'}\to\Delta^m$ takes one simplex to the
other.

To see this, we repeat the steps of the construction. The first
step was the definition of the map $\phi^*:A\to\bar A$. The
definition is compatible with structure maps because the
bimorphism $\phi:u\to\bar u$ is uniquely defined.

The second step was the choice of pushout diagrams~\eqref{push}.
But this choice is made in $\cc C=C^b(\cc A)$, and an element of
$S_m\cc C$ is a certain diagram in $\cc C$ on which the simplicial
structure maps in $S.\cc C$ operate by omission and/or
reduplicating of data. So again there is the required
compatibility.

The third step was construction of the isomorphisms $(\bar
a_i,\bar c_i,\bar b_i):\bar E_{i-1}\to\bar E_i$. The desired
compatibility follows from the fact that the maps $(\bar
a_i,a_i,c_i)$ represented by diagrams~\eqref{33} are well defined.

With an extra care one can arrange the choices so that the
homotopy starts from the identity map (namely if $A\to\bar A$ is
an identity map we insist that $C\to\bar C$ is also an identity
map); and that the image of $v_{n*}j_*$ is fixed under the
homotopy (namely if $\bar A=0$ we insist that $\bar C\to\bar B$ is
the identity map on $\bar B$). We have now constructed the desired
homotopy. This completes the proof.
\end{proof}

\renewcommand{\proofname}{Proof}

Now we discuss some immediate consequences of the additivity
theorem. Let $\cc A$ and $\cc A'$ be two exact categories. By an
{\it exact sequence\/} of exact functors $\cc A \to\cc A'$ is
meant a sequence of natural transformations $F'\to F\to F''$ such
that for every $A\in\cc A$ the sequence $F'(A)\to F(A)\to F''(A)$
is exact in $\cc A'$.

If $\cc A',\cc A''$ are fully exact subcategories of an exact
category $\cc A$ by $\cc E(\cc A',\cc A,\cc A'')$ denote the exact
subcategory of the exact extension category $\cc E=\cc E(\cc A)$
with the source and target entries in $\cc A'$ and $\cc A''$
respectively.

\begin{prop}[Equivalent formulations of the additivity theorem]\label{eqq}

Each of the following conditions is equivalent to the additivity
theorem (Theorem~\ref{add}).

\begin{enumerate}
\item The following projection is a homotopy equivalence,
   $$i.\bb S.\cc E(\cc A',\cc A,\cc A'')\to i.\bb S.\cc A'\times i.\bb S.\cc A'',
     \ \ \ A\rightarrowtail C\twoheadrightarrow B\longmapsto (A,B).$$

\item The following two maps are homotopic,
   $$i.\bb S.\cc E\to i.\bb S.\cc A,
     \ \ \ A\rightarrowtail C\twoheadrightarrow B\longmapsto C,\textrm{ respectively }
     A\ps B.$$

\item If $F'\to F\to F''$ is an exact sequence of exact functors
$\cc A\to\cc A'$ then there exists a homotopy
   $$|i.\bb S.F|\iso |i.\bb S.F'|\vee|i.\bb S.F''|.$$
\end{enumerate}
\end{prop}

\begin{proof}
The proof is similar to that of~\cite[1.3.2]{W}.
\end{proof}

If $F:\cc A\to\cc C$ and $G:\cc B\to\cc C$ be two arbitrary
functors with common codomain, the {\it fibre product\/}
$\prod(F,G)$ is defined as the category of triples
   $$(A,c,B),A\in\cc A,B\in\cc B,c:F(A)\to G(B)\textrm{ is an isomorphism.}$$
This is equivalent to the pullback category in special cases, for
example if one of $F$ and $G$ is a retraction, but not in general.
It follows from~\cite[p.~180]{W1} that if $F$ and $G$ are exact
functors then $\prod(F,G)$ is an exact category in a natural way,
and the projections to $\cc A$ and $\cc B$ are exact functors.
This is directly extended to the definition of the fibred product
for simplicial exact functors of simplicial exact categories.

Let $\cc B$ be an exact category and $P(S.\cc B)$ be the path
space for the simplicial exact category $S.\cc B$. One has the
simplicial map $\partial_0:P(S.\cc B)\to S.\cc B$. Let $F:\cc
A\to\cc B$ be an exact functor with $\cc A$ an exact category.
Denote by $S.(\cc A\to\cc B)$ the fibred product of the diagram
   $$S.\cc A\lra{F}S.\cc B\bl{\partial_0}\longleftarrow P(S.\cc B).$$
Each $S_n(\cc A\to\cc B)$, $n\geq 0$, consists of the triples
   $$(A,c,B),A\in S_n\cc A,B\in S_{n+1}\cc B,c:FA\lra{\sim}\partial_0B.$$
It is an exact category by above.

Considering $\cc B$ as a simplicial category in a trivial way we
have an inclusion $\cc B\to P(S.\cc B)$ whose composition with the
projection to $S.\cc B$ is trivial. Lifting this inclusion to
$S.(\cc A\to\cc B)$, and combining with the other projection, we
then obtain a sequence
   $$D^b(\cc B)\to\bb S.(\cc A\to\cc B)\to\bb S.\cc A$$
in which $\bb S.(\cc A\to\cc B)=D^b(S.(\cc A\to\cc B))$ and the
composed map is trivial.

\begin{prop}\label{j}
The sequence
   $$i.\bb S.\cc B\lra{}i.\bb S.\bb S.(\cc A\to\cc B)\lra{}i.\bb S.\bb S.\cc A,$$
in which $\bb S.\bb S.(\cc A\to\cc B)=D^b(S.(S.\cc A\to S.\cc B))$
is a fibration up to homotopy.
\end{prop}

\begin{proof}
The proof is similar to that of~\cite[1.5.5]{W}.
\end{proof}

Similarly, there is a sequence
   $$i.\bb S.\cc B\to P(i.\bb S.\bb S.\cc B)\to i.\bb S.\bb S.\cc B$$
where the ``$P$" refers to the first $S.$-direction, say. As a
special case of the preceding proposition this sequence is a
fibration up to homotopy.

Thus $|i.\bb S.\cc B|\to\Omega|i.\bb S.\bb S.\cc B|$ is a homotopy
equivalence and more generally, in view of Lemma~\ref{segal}, also
the map $|i.\bb S.^n\cc B|\to\Omega|i.\bb S.^{n+1}\cc B|$ for
every $n\geq 1$, proving the postponed claim that the spectrum
$n\longmapsto|i.\bb S.^n\cc B|$ is a $\Omega$-spectrum beyond the
first term.

\begin{cor}\label{111}
Suppose we are given a sequence $\cc A\to\cc B\to\cc C$ of exact
functors between exact categories. Then the square
   $$\begin{CD}
      i.\bb S.\cc B@>>>i.\bb S.\bb S.(\cc A\to\cc B)\\
      @VVV@VVV\\
      i.\bb S.\cc C@>>>i.\bb S.\bb S.(\cc A\to\cc C)
     \end{CD}$$
is homotopy cartesian.
\end{cor}

\begin{proof}
See~\cite[1.5.6]{W}.
\end{proof}

\begin{cor}\label{333}
The following two assertions are valid.

$(1)$ To an exact functor there is associated a fibration
   $$i.\bb S.\cc B\to i.\bb S.\cc C\to i.\bb S.\bb S.(\cc B\to\cc C).$$

$(2)$ If $\cc C$ is a retract of $\cc B$ (by exact functors) there
is a splitting
   $$i.\bb S.\cc B\iso i.\bb S.\cc C\times i.\bb S.\bb S.(\cc C\to\cc B).$$
\end{cor}

\begin{proof}
\cite[1.5.7,~1.5.8]{W}.
\end{proof}

By a {\it nice complicial biWaldhausen category\/} $\cc C$ formed
from the category of complexes $C^b(\cc B)$ with $\cc B$ an
abelian category will be meant a complicial biWaldhausen category
in the sense of Thomason~\cite{T} which is closed under the
formation of canonical homotopy pushouts and canonical homotopy
pullbacks. For example, let $\cc A$ be an exact category and $\cc
A\to\cc B$ be the Gabriel-Quillen embedding~\cite[Appendix~A]{T}.
Then $C^b(\cc A)$ is a nice complicial biWaldhausen category
formed from $C^b(\cc B)$.

It is directly verified that for any $n$ the category $S_n\cc C$
is a nice complicial biWaldhausen category which is formed from
the category of complexes $C^b(\cc B^{\Ar\Delta^n})$ with $\cc
B^{\Ar\Delta^n}$ the abelian functor category
$\Hom(\Ar\Delta^n,\cc B)$. The relevant subcategories of
bifibrations and weak equivalences are defined componentwise. In a
similar way, given a small category $I$ the diagram category $\cc
C^I$ is a nice complicial biWaldhausen category. There results a
simplicial nice complicial biWaldhausen category
   $$S.:\Delta^n\longmapsto S_n\cc C.$$

Let $w^{-1}\cc C$ denote the derived category obtained from $\cc
C$ by inverting weak equivalences. It is canonically triangulated
and the homotopy category $\cc C/\iso$ admits both left and right
calculus of fractions~\cite[p.~269]{T}. One obtains the following
bisimplicial object
   $$i.\bb S.:\Delta^m\times\Delta^n\longmapsto i_mw^{-1}S_n\cc C.$$
Denote by $\cc E(\cc C)$ the extension category of $\cc C$. Then
the proof of the following statement is similar to complexes (all
tricks of paragraph~\ref{rid} are also valid for this case).

\begin{cor}\label{comp}
Let $\cc C$ be a nice complicial biWaldhausen category. Then the
map
   $$i.\bb S.\cc E(\cc C)\xrightarrow{(s_*,q_*)}i.\bb S.\cc C\times i.\bb S.\cc C$$
is a homotopy equivalence.
\end{cor}

\section{D\'erivateurs associated to complicial biWaldhausen categories}\label{c}

In this section we show that the additivity theorem is valid for
d\'erivateurs associated to nice complicial biWaldhausen
categories. Logically, one should now read Addendum, and then
return to this section.

Let $\cc C$ be a nice complicial biWaldhausen category. One of the
most important for applications d\'erivateurs (of the domain
$\dirf$) is given by the hyperfunctor
   $$\bb D\cc C:I\in\dirf\longmapsto w^{-1}\cc C^I$$
with $w^{-1}\cc C^I$ the derived category of the nice complicial
biWaldhausen diagram category $\cc C^I$~\cite{C1,K}. If $\cc
C=C^b(\cc A)$ with $\cc A$ an exact category, the corresponding
d\'erivateur is denoted by $\bb D^b(\cc A)$.

\begin{defs}{\rm
A left pointed d\'erivateur $\bb D$ of the domain $\ord$ is {\it
represented\/} by a nice complicial biWaldhausen category $\cc C$
if there is a right exact equivalence $F:\bb D\cc C\lra{}\bb D$.
This equivalence induces a homotopy equivalence of bisimplicial
sets $F:i.S.\bb D\cc C\lra{}i.S.\bb D$. If $\cc C=C^b(\cc A)$ with
$\cc A$ an exact category, we shall say that $\bb D$ is {\it
represented\/} by $\cc A$. }\end{defs}

\begin{lem}\label{fn}
The inclusion $F_n\cc C\to \cc C^{\Delta^n}$ induces an
equivalence of derived categories $w^{-1}F_n\cc
C\lra{\sim}w^{-1}\cc C^{\Delta^n}$.
\end{lem}

\begin{proof}
Since a map in $F_n\cc C$ is a weak equivalence if and only if it
is so in $\cc C^{\Delta^n}$, it suffices to show that given an
object $A=A_0\xrightarrow{}\cdots\xrightarrow{}A_n$ in $\cc
C^{\Delta^n}$ there is a quasi-isomorphism from an object in
$F_n\cc C$ to $A$. Let us consider the following diagram in $\cc
C$.

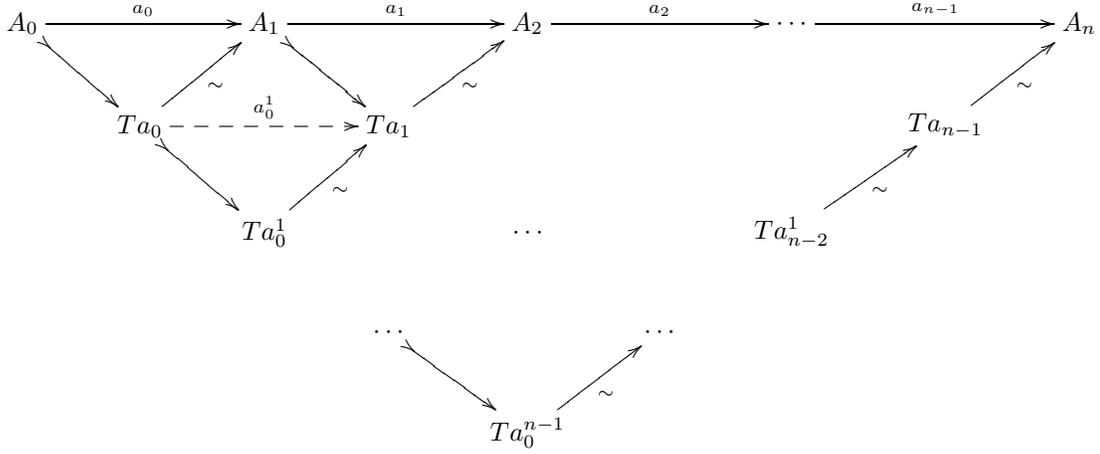
\footnotesize

   $$\xymatrix{
     A_0\ar[rr]^{a_0}\ar@{ >->}[dr]&& A_1\ar[rr]^{a_1}\ar@{ >->}[dr]&& A_2
     \ar[rr]^{a_2}&&\cdots\ar[rr]^{a_{n-1}}&&A_n\\
     &Ta_0\ar@{ >->}[dr]\ar@{-->}[rr]^{a^1_0}\ar[ur]_{\sim}&&Ta_1\ar[ur]_{\sim}&&&&Ta_{n-1}\ar[ur]_{\sim}\\
     &&Ta_0^1\ar[ur]_{\sim}&&\cdots&&Ta^1_{n-2}\ar[ur]_{\sim}\\
     &&&\cdots\ar@{ >->}[dr]&&\cdots\\
     &&&&Ta^{n-1}_0\ar[ur]_{\sim}}$$
\normalsize Here $T(-)$ stands for the cylinder object of a
morphism, the arrows labelled with ``$\sim$" are weak
equivalences, and all the squares of the diagram are commutative.
This diagram yields a weak equivalence from the object
$A_0\rightarrowtail Ta_0\rightarrowtail\cdots\rightarrowtail
Ta^{n-1}_0$ of $F_n\cc C$ to the object
$A_0\xrightarrow{}\cdots\xrightarrow{}A_n$ of $\cc C^{\Delta^n}$,
whence the assertion.
\end{proof}

\begin{lem}\label{zu}
If $\bb D$ is a d\'erivateur represented by a nice complicial
biWaldhausen category $\cc C$, then $\bb S_n\bb D$ is represented
by the nice complicial biWaldhausen category $S_n\cc C$ for all
$n$.
\end{lem}

\begin{proof}
The image of each cocartesian square in $\cc C$ with two parallel
arrows cofibrations is cocartesian in $\bb D\cc C_\Box$. This
yields a right exact morphism $\bb DS_n\cc C\to\bb S_n\bb D\cc C$.
Consider the commutative diagram of left pointed d\'erivateurs
   $$\begin{CD}
      \bb DF_{n-1}\cc C@>>>\bb D\cc C(\Delta^{n-1})\\
      @A\sim AA@AA\sim A\\
      \bb DS_n\cc C@>>>\bb S_n\bb D\cc C@>{\sim}>>\bb S_n\bb D
     \end{CD}$$
in which the morphisms marked with ``$\sim$" are right exact
equivalences. The equivalence on the left is induced by the
equivalence $S_n\cc C\to F_{n-1}\cc C$ and the equivalence on the
right is a consequence of~\cite[3.1]{Gar}.

The upper arrow is an equivalence by Lemma~\ref{fn}. We see from
the commutative diagram above that $\bb DS_n\cc C\to\bb S_n\bb
D\cc C$ is a right exact equivalence. Therefore the composed map
$\bb DS_n\cc C\to\bb S_n\bb D$ produces a Waldhausen model for
$\bb S_n\bb D$.
\end{proof}

\begin{cor}\label{trr}
A natural map of bisimplicial sets $i.\bb S.\cc C\to i.S.\bb D\cc
C$ is a homotopy equivalence.
\end{cor}

\begin{proof}
By the proof of the preceding lemma the functor $w^{-1}S_n\cc C\to
S_n\bb D\cc C$ is an equivalence of categories for every $n$.
Lemma~\ref{segal} implies the claim.
\end{proof}

Let $\cc E$ denote the extension category of $\cc C$.

\begin{cor}\label{cvv}
If $\bb D$ has a Waldhausen model then so does $\bb E=\bb E(\bb
D)$.
\end{cor}

\begin{proof}
It is enough to consider the commutative diagram
   $$\begin{CD}
      \bb DS_2\cc C@>\sim>>\bb S_2\bb D\\
      @V\sim VV@VV\sim V\\
      \bb D\cc E@>>>\bb E
     \end{CD}$$
in which the left arrow is a right exact equivalence by the exact
equivalence of $S_2\cc C$ and $\cc E$, the upper arrow is a right
exact equivalence by the preceding lemma, and the right arrow is a
right exact equivalence by~\cite[6.2]{Gar}.
\end{proof}

There are three natural right exact maps $s,t,q:\bb E\to\bb D$
taking $E$ to $E_{(0,0)}$, $E_{(0,1)}$ and $E_{(1,1)}$
respectively. The following result states that the additivity
theorem holds for d\'erivateurs represented by nice complicial
biWaldhausen categories.

\begin{thm}\label{mod}
Let $\bb D$ be a d\'erivateur represented by a nice complicial
biWaldhausen category. Then the map
   $$i.S.\bb E\xrightarrow{(s_*,q_*)}i.S.\bb D\times i.S.\bb D$$
is a homotopy equivalence.
\end{thm}

\begin{proof}
Let $\bb D\cc C\lra{\sim}\bb D$ be a Waldhausen model for $\bb D$
and let $\cc E$ be the extension category of $\cc C$. By
Corollary~\ref{cvv} $\bb D\cc E\lra{\sim}\bb E$ is a Waldhausen
model for $\bb E$. Consider the following commutative diagram
   $$\begin{CD}
      i.\bb S.\cc E@>(s_*,q_*)>>i.\bb S.\cc C\times i.\bb S.\cc C\\
      @VVV@VVV\\
      i.S.\bb D\cc E@>(s_*,q_*)>>i.S.\bb D\cc C\times i.S.\bb D\cc C\\
      @VVV@VVV\\
      i.S.\bb E@>(s_*,q_*)>>i.S.\bb D\times i.S.\bb D.
     \end{CD}$$
The vertical arrows are homotopy equivalences. By
Corollary~\ref{comp} the map $i.\bb S.\cc E\to i.\bb S.\cc C\times
i.\bb S.\cc C$ is a homotopy equivalence. This implies the claim.
\end{proof}

Let $\bb D$ be a d\'erivateur represented by a nice complicial
biWaldhausen category. We can apply the $S.$-construction to each
$\bb S_n\bb D$, obtaining a bisimplicial left pointed d\'erivateur
represented by a nice complicial biWaldhausen category. Iterating
this construction, we can form the multisimplicial object $\bb
S.^n\bb D=\bb S.\bb S.\cdots\bb S.\bb D$ and the multisimplicial
categories $iS.^n\bb D$ of isomorphisms. The assertion below shows
that $|i.S.^n\bb D|$ is the loop space of $|i.S.^{n+1}\bb D|$ for
any $n\geq 1$ and that the sequence
   $$\Omega|i.S.\bb D|,\Omega|i.S.S.\bb D|,\ldots,\Omega|i.S.^n\bb D|,\ldots$$
forms a connective $\Omega$-spectrum $\bb K\bb D$ (see the
definition of the structure maps in~\cite{Gar}). In this case, one
can think of the $K$-theory of $\bb D$ in terms of this spectrum.
This does not affect the $K$-groups, because:
   $$\pi_i(\bb K\bb D)=\pi_i(K(\bb D))=K_i(\bb D),\ \ \ i\geq 0.$$

\begin{cor}\label{spektr}
Let $\bb D$ be a d\'erivateur represented by a nice complicial
biWaldhausen category. Then
   $$n\longmapsto i.S.^n\bb D$$
is a $\Omega$-spectrum beyond the first term. In particular, the
$K$-theory of $\bb D$ can equivalently be defined as the space
   $$\Omega^\infty|i.S.^\infty\bb D|=\lim_n\Omega^n|i.S.^n\bb D|.$$
\end{cor}

\begin{proof}
For every $n\geq 0$, it follows from Lemma~\ref{zu} that $\bb
S_n\bb D$ is a d\'erivateur represented by a nice complicial
biWaldhausen category. By Theorem~\ref{mod} the class of such
d\'erivateurs satisfies the addivity theorem. The claim now
follows from~\cite[section~6]{Gar}.
\end{proof}

Let us say a few words what thing goes wrong when conforming
Waldhausen's~\cite{W} or McCarthy's~\cite{Mc} proof of additivity
--- in fact, both have the same complexity --- to d\'erivateurs
(or systems of diagram categories).

The first step is to show that additivity follows from the fact
that for any Waldhausen category $\cc C$ the map
   $$s.\cc E(\cc C)\to s.\cc C\times s.\cc C$$
with $s.\cc C=\Ob S.\cc C$ is a homotopy equivalence (just apply
the same map to the Waldhausen category $w_n\cc C, n\geq 0,$ of
strings of weak equivalences and then apply Lemma~\ref{segal}).
The same applies to left pointed d\'erivateurs: it suffices to
show that for any left pointed d\'erivateur $\bb D$ the map
   $$s.\bb E(\bb D)\to s.\bb D\times s.\bb D$$
with $s.\bb D=\Ob S.\bb D$ is a homotopy equivalence. This is
because the hyperfunctor $I\longmapsto i_n\bb D_I$ taking an index
category to the category of strings of isomorphisms is a left
pointed d\'erivateur, too.

The second step consists of choices of pushouts (we neglect
quotients here). Precisely, we are given two maps with common
source $f:A\to C$ and $\phi:A\to\bar A$ and $f$ a cofibration
representing some simplex. Afterwards one constructs a pushout
square
   $$\begin{CD}
      A@>f>>C\\
      @V{\phi}VV@VVV\\
      \bar A@>{\bar f}>>\bar C
     \end{CD}$$
to get a simplex $\bar A\rightarrowtail\bar C$ from the simplex
$A\rightarrowtail C$. However it is not immediately clear that the
same procedure applies to d\'erivateurs. Normally we are given, as
above, two objects $X,Y$ in $\bb D_{\Delta^1}$ with common
``source" $X_0=Y_0$ and $X$ some simplex. To get a ``new" simplex
$\bar X$ from $X$ and $Y$ in the same way one should construct a
cocartesian square whose projection on $(0,0)\to(0,1)$ is the $X$
and that on $(1,0)\to(1,1)$ is the $\bar X$. It seems that
d\'erivateurs do not have enough information to do so. We were
able to do that for d\'erivateurs represented by nice complicial
biWaldhausen categories by using certain tricks of
paragraph~\ref{rid}, but it is not clear how to construct the
necessary homotopy for all (left pointed) d\'erivateurs basing
only on the known proofs of additivity for Waldhausen categories.

To conclude this section, we would like to invite experts to prove
additivity for d\'erivateurs represented by closed model
categories. A similar technique used in this paper should be
applicable for this case as well.

\section{The derived $K$-theory of an exact category}\label{5}

In this section we define the derived $K$-theory $DK(\cc A)$ of an
exact category $\cc A$. Though it is homotopy equivalent to the
$K$-theory of its d\'erivateur $\bb D^b(\cc A)$ it is more
convenient for practical reasons to deal with the space $DK(\cc
A)$ than with the space $K(\bb D^b(\cc A))$.

\begin{defs}{\rm
The {\it Algebraic $DK$-theory\/} of an exact category $\cc A$ is
defined as the pointed space
   $$DK(\cc A)=\Omega|i.\bb S.\cc A|.$$
The $DK$-groups of $\cc A$ are the homotopy groups of $DK(\cc A)$
   $$DK_*(\cc A)=\pi_*(\Omega|i.\bb S.\cc A|)=\pi_{*+1}(|i.\bb S.\cc A|).$$
}\end{defs}

Let $(ExCats)$ denote the category of exact categories and exact
functors. There results a functor
   $$DK:(ExCats)\to (Spaces)$$

It follows from Corollary~\ref{trr} that a natural map $DK(\cc
A)\to K(\bb D^b(\cc A))$ is a homotopy equivalence. Hence the
$DK$-theories $DK(\cc A)$ and $DK(\cc A')$ of exact categories
$\cc A$ and $\cc A'$ are homotopy equivalent whenever $K(\bb
D^b(\cc A))$ and $K(\bb D^b(\cc A'))$ are. It also follows
(see~\cite{Gar}) that $DK_0(\cc A)$ is isomorphic to the
Grothendieck group $K_0(\cc A)$.

We prove below some basic results about $DK$-theory. The first
result is the Additivity Theorem.

\begin{thm}[Additivity]\label{123}
Let $\cc A$ be an exact category and $\cc E$ its extension
category. Then the map
   $$DK(s,q):DK(\cc E)\to DK(\cc A)\times DK(\cc A)$$
is a homotopy equivalence. If $F'\to F\to F''$ is an exact
sequence of exact functors $\cc A\to\cc A'$ then there is a
homotopy of maps
   $$DK(F)\iso DK(F')\vee DK(F''):DK(\cc A)\to DK(\cc A').$$
The $DK$-theory of $\cc A$ can equivalently be defined as the
space
   $$\Omega^\infty|i.\bb S.^\infty\cc A|= \lim_n\Omega^n|i.\bb S.^n\cc A|.$$
One can also think of the $DK$-theory in terms of the
$\Omega$-spectrum
   $$\Omega|i.\bb S.\cc A|,\Omega|i.\bb S.\bb S.\cc A|,\ldots,\Omega|i.\bb S.^n\cc A|,\ldots$$
\end{thm}

\begin{proof}
These follow from results of section~\ref{234}.
\end{proof}

\subsection{Approximation and Resolution theorems}

In this paragraph we prove a modified Approximation Theorem and
Resolution Theorem.

\begin{thm}[Approximation]\label{app}
Let $\cc A$ and $\cc A'$ be two exact categories and let $w\cc C$
and $w\cc C'$ denote the Waldhausen categories of
quasi-isomorphisms in $\cc C=C^b(\cc A)$ and in $\cc C'=C^b(\cc
A')$ respectively. Suppose further that $F:w\cc C\to w\cc C'$ is
an exact functor of Waldhausen categories such that it induces an
equivalence of the derived categories $D^b(\cc A)\lra{\sim}D^b(\cc
A')$. Then $DK(\cc A)$ is homotopy equivalent to $DK(\cc A')$. If
$F$ is induced by an exact functor $f:\cc A\to\cc A'$, this
homotopy equivalence is given by the induced map $DK(f):DK(\cc
A)\to DK(\cc A')$.
\end{thm}

We postpone the proof and first define some new concepts and prove
certain technical lemmas.

\begin{defs}{\rm

Under the notation of Theorem~\ref{app} we say that $F$ has the
{\it approximation property\/} (respectively {\it h-approximation
property}) if it meets the axioms $App1$ and $App2$ below
(respectively the axioms $App1$ and $HApp2$).

\begin{itemize}

\item[$App1$] A map in $\cc C$ is a quasi-isomorphism if and only
if its image is a quasi-isomorphism in $\cc C'$.

\item[$App2$] Any map $f:FX\to Y$ in $\cc C'$, $X\in\cc C$, fits
into a commutative diagram
   $$\begin{CD}
      FX@>f>>Y\\
      @VFpVV@VVsV\\
      FX'@>t>>Y'
     \end{CD}$$
in which $p:X\to X'$ is a map in $\cc C$ and $s,t$ are
quasi-isomorphisms in $\cc C'$.

\item[$HApp2$] Any map $f:FX\to Y$ in $\cc C'$, $X\in\cc C$, fits
into a homotopy commutative diagram
   $$\begin{CD}
      FX@>f>>Y\\
      @VFpVV@VVsV\\
      FX'@>t>>Y'
     \end{CD}$$
with $s$ and $t$ quasi-isomorphisms in $\cc C'$.

\end{itemize}

$F$ has the {\it dual approximation property\/} (respectively {\it
dual h-appro\-xi\-ma\-tion property}) if the axiom $App1$ and the
dual axiom $App2^{\op}$ obtained by reversing the direction of
arrows in $App2$ (respectively the axiom $App1$ and the dual axiom
$HApp2^{\op}$) are satisfied. These are a modification for the
Waldhausen axioms $WApp1-WApp2$~\cite[p.~352]{W}.

}\end{defs}

The next statement is due to Cisinski~\cite{C1}. In fact, he
proves it in a more general setting.

\begin{prop}\label{spi}
Under the notation of Theorem~\ref{app} the following are
equivalent:

$(1)$ the functor $F$ induces an equivalence of the derived
categories $D^b(\cc A)\to D^b(\cc A')$;

$(2)$ the functor $F$ has the approximation property;

$(3)$ the functor $F$ has the h-approximation property;

$(2^{\op})$ the functor $F$ has the dual approximation property;

$(3^{\op})$ the functor $F$ has the dual h-approximation property.
\end{prop}

\begin{proof}
It is enough to show
$(1)\Longleftrightarrow(2^{\op})\Longleftrightarrow(3^{\op})$ (the
equivalences $(1)\Longleftrightarrow(2)\Longleftrightarrow(3)$ are
dually proved).

$(2^{\op})\Longrightarrow(3^{\op})$ is obvious. Let us show
$(3^{\op})\Longrightarrow(2^{\op})$. Any arrow $f:Y\to FX$ can be
fitted into the following diagram:
   $$\xymatrix{&& FX''\ar@{.>}[dl]^{Fl}\ar@/^/@{.>}[ddl]^{Fg'}\\
               \ar@{}[dr] |{\textrm{h. comm.}}
               Y' \ar@/^/@{.>}[urr]^{t'}\ar[d]_{s} \ar[r]^{t} & FX' \ar[d]^{Fg}\\
               Y \ar[r]_{f} & FX}$$
where the homotopy commutative square with the entries
$(Y',FX',Y,FX)$ exists by assumption, $X''=\Cocyl(g)$ and the
square with the entries $(Y',FX'',Y,FX)$ is commutative (see
paragraph~\ref{rid}).

$(1)\Longrightarrow(3^{\op})$. The axiom $App1$ is obvious. Check
the axiom $HApp2$. Let $\alpha:Y\to FX$ be an arrow in $\cc C'$.
There exists an isomorphism $ts^{-1}:FX'\to Y$ in $D^b(\cc A')$
resulting a map $FX'\to FX$. Let this map be the image of a map
$X'\to X$ in $D^b(\cc A)$ represented by a diagram $X'\bl
q\longleftarrow X''\lra{f} X$.

There is a common denominator
   $$\xymatrix{&&Y''\ar[dl]_{s'}\ar[dr]^{f'}\\
      &U\ar[dl]_{s}\ar[drrr]^(.2){\alpha t}&& FX''\ar[dlll]_(.2){Fq}\ar[dr]^{Ff}\\
      FX'&&&&FX}$$
yielding a homotopy commutative square
   $$\xymatrix{\ar @{}[dr] |{\textrm{h. comm.}}
               Y'' \ar[d]_{ts'} \ar[r]^{f'} & FX'' \ar[d]^{Ff} \\
               Y \ar[r]_{\alpha}      & FX}$$

It remains thus to show $(2^{\op})\Longrightarrow(1)$. Given an
object $Y\in D^b(\cc A')$ there is a diagram
   $$\xymatrix{Y'\ar[d]_{s} \ar[r]^{t} & FX \ar[d]\\
               Y \ar[r] & 0}$$
with $s,t$ quasi-isomorphisms. We see that $D^b(\cc A)\to D^b(\cc
A')$ is essentially surjective. Let us show that it is full.

Let $\alpha:FX\to FX'$ be a map in $D^b(\cc A')$ represented by a
diagram $FX\bl s\longleftarrow Y\lra{f}FX'$. There is a
commutative diagram
   $$\begin{CD}
      Y'@>q>>FZ\\
      @VtVV@VV{(Fu,Fv)^t}V\\
      Y@>{(s,f)^t}>>F(X\times X')
     \end{CD}$$
with $q,t$ quasi-isomorphisms. It follows that $Fu$ is a
quasi-isomorphism. We get a denominator
   $$\xymatrix{&&Y'\ar[dl]_{t}\ar[dr]^{q}\\
      &Y\ar[dl]_{s}\ar[drrr]^(.2){f}&& FZ\ar[dlll]_(.2){Fu}\ar[dr]^{Fv}\\
      FX&&&&FX'}$$
This shows that $\alpha$ is the image of $vu^{-1}:X\to X'$.

To show that the functor in question is faithful we shall need the
following

\begin{sublem}
Let two maps $u,v:X\to Y$ in $\cc C$ be such that there is a
quasi-isomorphism $q:U\to FX$ and $Fu\circ q$ and $Fv\circ q$ are
homotopic in $\cc C'$. Then there exists a quasi-isomorphism
$s:T'\to X$ in $\cc C$ such that $us$ is homotopic to $vs$.
\end{sublem}

\begin{proof}
Let $\Cocyl(Y)$ denote the cocylinder of the map $1_Y$ and write
$d_0,d_1:\Cocyl(Y)\to Y$ for the natural projections. The map
$(d_0,d_1):\Cocyl(Y)\to Y\times Y$ is an epimorphism in $\cc C$.
Construct a cartesian diagram
   $$\begin{CD}
      T@>{\sigma}>>\Cocyl(Y)\\
      @V{\pi}VV@VV{(d_0,d_1)}V\\
      X@>{(u,v)}>>Y\times Y.
     \end{CD}$$
Since $Fu\circ q,Fv\circ q$ are homotopic, there is a map
$\alpha:U\to F[\Cocyl(Y)]=\Cocyl(FY)$ such that
$(Fd_0,Fd_1)\circ\alpha=(Fu,Fv)\circ q$. There results a
commutative diagram
   $$\begin{CD}
      U@>{\alpha}>>F[\Cocyl(Y)]\\
      @V{q}VV@VV{(Fd_0,Fd_1)}V\\
      FX@>{(Fu,Fv)}>>FY.
     \end{CD}$$
There is a unique map $v:U\to FT$ such that $F(\sigma)\circ
v=\alpha$ and $F(\pi)\circ v=q$.

The map $v$ fits into a commutative square
   $$\begin{CD}
      U'@>{c}>>FT'\\
      @V{d}VV@VV{F\delta}V\\
      U@>{v}>>FT
     \end{CD}$$
with $c,d$ quasi-isomorphisms. We see that $F(\pi\delta)\circ
c=qd$ is a quasi-isomorphism, and so is $F(\pi\delta)$. By $App1$
the map $\pi\delta$ is a quasi-isomorphism. It follows that
$u(\pi\delta)$ is homotopic to $v(\pi\delta)$, hence the required
map $s$ is $\pi\delta$.
\end{proof}

Now prove that the functor $D^b(\cc A)\to D^b(\cc A')$ is
faithful. It suffices to show that if $Ff=0$ in $D^b(\cc A')$ with
$f:X\to Y$ a map in $\cc C$, then $f$ equals to zero in $D^b(\cc
A)$.

The property for a map $\alpha$ in $D^b(\cc A')$ of being equal to
zero is equivalent to saying that there is a quasi-isomorphism $q$
such that $\alpha q$ is homotopic to zero. By our assumption
$Ff\circ q\sim 0$, and so there exists a quasi-isomorphism $s$ in
$\cc C$ with $fs\sim 0$ by the sublemma above. It follows that $f$
equals to zero in $D^b(\cc A)$. We are done.
\end{proof}

The last proposition also applies to nice complicial biWaldhausen
categories $\cc A$ and $\cc B$ and a complicial exact functor
between them. The following shows to which extent the Thomason
Approximation Theorem~\cite[1.9.8]{T} for nice complicial
biWaldhausen categories is a modification of the Waldhausen
Approximation Theorem~\cite[1.6.7]{W}

\begin{cor}\label{tho}
Let $\cc A$ and $\cc B$ be nice complicial biWaldhausen categories
and let $F:\cc A\to\cc B$ be a complicial exact functor. Suppose
that $F$ has the approximation or the h-approximation property
(the axioms $App1-App2$ or $App1-HApp2$). Then $F$ induces a
homotopy equivalence of $K$-theory spaces
   $$K(F):K(\cc A)\to K(\cc B).$$
\end{cor}

\begin{proof}
By Proposition~\ref{spi} $F$ induces an equivalence of the derived
categories $w^{-1}F:w^{-1}\cc A\to w^{-1}\cc B$. Then $K(F)$ is a
homotopy equivalence by~\cite[1.9.8]{T}.
\end{proof}

\begin{lem}\label{nuu}
An exact functor $F:w\cc C\to w\cc C'$ meets the axiom $App2$ if
and only if any arrow $f:FX\to Y$ in $\cc C'$, $X\in\cc C$, fits
into a commutative diagram
   $$\xymatrix{FX\ar@{ >->}[d]_{Fj_1}\ar[r] & Y \ar@{ >->}[d]^t\\
               FX'\ar[r]_s & Y'}$$
in which $j_1,t$ are cofibrations and $s,t$ are quasi-isomorphisms
in $\cc C'$.
\end{lem}

\begin{proof}
If $F$ meets the axiom $App2$, then any arrow $f:FX\to Y$ in $\cc
C'$ fits into a commutative diagram
   $$\begin{CD}
      FX@>f>>Y\\
      @VFaVV@VVpV\\
      FX'@>q>>Y'
     \end{CD}$$
in which $p,q$ are quasi-isomorphisms. Let $T=\Cyl(a)$, then
$a=rj_1$ with $j_1:X\to T$ a cofibration and $r:T\to X'$ a
quasi-isomorphism. Construct a cocartesian square
   $$\xymatrix{FX\ar@{ >->}[d]_{Fj_1}\ar[r]^f & Y \ar@{ >->}[d]^{u}\\
               FT\ar[r]_{f'} & V}$$
There exists a unique map $v:V\to Y'$ such that $q\circ Fr=vf'$
and $vu=p$. The map $v$ factors as $V\lra{m}T'\lra{n}Y'$ with
$T'=\Cyl(v)$, $m$ a cofibration and $n$ a quasi-isomorphism. There
results a commutative diagram
   $$\xymatrix{FX\ar@{ >->}[r]^{Fj_1}\ar[d]_f & FT\ar[d]^{mf'}\ar[r]^{Fr}& FX'\ar[d]^q\\
               Y \ar@{ >->}[r]^{mu} & V\ar[r]^n& Y'}$$
Since $p=nmu$ and $p,n$ are quasi-isomorphisms, we see that $mu$
is a quasi-iso\-mor\-phism. Also, $mf'$ is a quasi-isomorphism
because $n,q,Fr$ are. We are done.
\end{proof}

\renewcommand{\proofname}{Proof of Theorem~\ref{app}}

\begin{proof}
In view of Lemma~\ref{segal} and Proposition~\ref{app} it suffices
to check that for any $n\geq 0$ the induced functor $wS_n\cc C\to
wS_n\cc C'$ has the approximation property. This is so for $n=0$.
Obviously, it is enough to check the approximation property for
the map $wF_n\cc C\to wF_n\cc C'$ and $n\geq 0$. For $n=0$ it
follows from our assumption and Proposition~\ref{app}. If we show
this for $n=1$ the general case will follow by induction.

Let the diagram represent a map $a:FX\to Y$ in $F_1\cc
C=C^b(F_1\cc A)$,
   $$\xymatrix{FX_0\ar[d]_{a_0} \ar@{ >->}[r] & FX_1 \ar[d]\\
               Y_0 \ar@{ >->}[r] & Y_1}$$
Then $a_0$ fits into a commutative square
   $$\xymatrix{FX_0\ar[d]_{Fq_0} \ar[r]^{a_0} & Y_0 \ar[d]\\
               FX'_0 \ar[r] & Y'_0}$$
One can construct the following commutative diagram
   $$\xymatrix@!0{
     &Y_0\ar@{ >->}[rr]\ar'[d][dd] && Y_1\ar[dd]\\
     FX_0\ar[ur]^{a_0}\ar@{ >->}[rr]\ar[dd]_{Fq_0} && FX_1\ar[ur]\ar[dd]\\
     &Y_0'\ar@{ >->}'[r][rr] && Y\\
     FX'_0\ar@{ >->}[rr]\ar[ur] && FX\ar[ur]_\alpha
     }$$
with $X=X_0'\coprod_{X_0}X_1$ and $Y=Y_0'\coprod_{Y_0}Y_1$. By
Lemma~\ref{nuu} $\alpha$ fits into a commutative diagram
   $$\xymatrix{FX\ar@{ >->}[d]_{Fj_1}\ar[r]^\alpha & Y \ar@{ >->}[d]^t\\
               FX'_1\ar[r]_s & Y'_1}$$
in which $j_1,t$ are cofibrations and $s,t$ are quasi-isomorphisms
in $\cc C'$. We get the commutative diagram
   $$\xymatrix@!0{
     &Y_0\ar@{ >->}[rr]\ar'[d][dd] && Y_1\ar[dd]\\
     FX_0\ar[ur]^{a_0}\ar@{ >->}[rr]\ar[dd]_{Fq_0} && FX_1\ar[ur]\ar[dd]\\
     &Y_0'\ar@{ >->}'[r][rr] && Y'_1\\
     FX'_0\ar@{ >->}[rr]\ar[ur] && FX_1'\ar[ur]_s
     }$$
that shows $App2$. The axiom $App1$ is obvious. The theorem is
proved.
\end{proof}

\renewcommand{\proofname}{Proof}

\begin{thm}[Resolution]\label{res}
Let $\cc P$ be an extension closed full exact subcategory of an
exact category $\cc M$. Assume further that

$(1)$ if $M'\rightarrowtail M\twoheadrightarrow M''$ is exact in
$\cc M$ and $M,M''\in\cc P$, then $M'\in\cc P$ and

$(2)$ for any object $M\in\cc M$ there is a finite resolution
$0\to P_n\to P_{n-1}\to\cdots\to P_0\to M\to 0$ with $P_i\in\cc
P$.

Then $DK(\cc P)\to DK(\cc M)$ is a homotopy equivalence (and thus
$DK_i(\cc P)\iso DK_i(\cc M)$ for all $i$).
\end{thm}

\begin{proof}
By~\cite[12.1]{K96} an extension closed full exact subcategory
$\cc P$ of an exact category $\cc M$ induces a fully faithful
canonical functor between their bounded derived categories if for
any exact sequence ${M''}\rightarrowtail{M'}\twoheadrightarrow P$
in $\cc M$ with $P\in\cc P$, there is a commutative diagram
   $$\begin{CD}
      0@>>>P''@>>>P'@>>>P@>>>0\\
      @.@VVV@VVV@|\\
      0@>>>M''@>>>M'@>>>P@>>>0
     \end{CD}$$
with $P',P''\in\cc P$ and in which the first row is also exact.
This condition follows from our assumptions. Indeed, by (2) one
can choose an admissible epimorphism $P'\to M'$ and the kernel of
the composed map $P'\to P$ is in $\cc P$ by (1).

Since each object $M\in\cc M$ has a finite resolution by objects
in $\cc P$, it follows that for every bounded complex $A$ with
entries in $\cc M$ there exists a quasi-isomorphism $B\to A$ for
some bounded complex $B$ with entries in $\cc P$ (the proof is
dual to that of~\cite[4.1(b)]{K90}). Therefore, the canonical
functor $D^b(\cc P)\to D^b(\cc M)$ is an equivalence.
Theorem~\ref{app} now implies the claim.
\end{proof}

\subsection{Pairings}

Let $\cc A,\cc B,\cc C$ be exact categories. We want to pair
Quillen's $K$-theory of $\cc A$ and the derived $K$-theory of $\cc
B$ into the derived $K$-theory of the latter. The appropriate
assumption to make is a pairing
   $$f:\cc A\times\cc B\to\cc C$$
which is a biexact functor, that is for each $A\in\cc A$ and
$B\in\cc B$ the partial functors
   $$f(A,-):\cc B\to\cc C,\ \ \ f(-,B):\cc A\to\cc C$$
are exact. We shall think of $f$ as a tensor product. For
technical reasons we assume that each of $\cc A,\cc B,\cc C$ is
equipped with a distinguished zero object 0 and that
$f(A,0)=0=f(0,B)$ always. Such a $f$ gives rise a pairing
   $$f:\cc A\times C^b(\cc B)\to C^b(\cc C)$$
which is also a biexact functor.

Let $s.\cc A$ denote the simplicial set $\{\Ob S_n\cc A\}_n$. We
obtain a map
   $$|s.\cc A|\times|i.\bb S.\cc B|\to|i.\bb S.\bb S.\cc C|$$
that takes $|s.\cc A|\vee|i.\bb S.\cc B|$ into the basepoint of
$|i.\bb S.\bb S.\cc C|$ because of the technical assumption we
made. This yields a map
   $$|s.\cc A|\wedge|i.\bb S.\cc B|\to|i.\bb S.\bb S.\cc C|$$
and hence a map
   $$\Omega|s.\cc A|\wedge\Omega|i.\bb S.\cc B|\to\Omega\Omega|i.\bb S.\bb S.\cc C|.$$
This is the desired pairing in $K$-theory in view of the homotopy
equivalence of $\Omega|s.\cc A|$ with $K(\cc A)$ and
$\Omega\Omega|i.\bb S.\bb S.\cc C|$ with $DK(\cc C)$. So we get a
map of spaces
   $$K(\cc A)\wedge DK(\cc B)\to DK(\cc C)$$
and hence homomorphisms of abelian groups
   $$K_m(\cc A)\otimes DK_n(\cc B)\to DK_{m+n}(\cc C),\ \ \ m,n\geq 0.$$

\subsection{Conjectures}

The central problem here is comparison of $DK(\cc A)$ with
Quillen's $K$-theory $K(\cc A)$. There is a natural map
$K(\rho):K(\cc A)\to DK(\cc A)$ factoring as
   $$K(\cc A)\xrightarrow{K(\tau)}K(wC^b(\cc
   A))\xrightarrow{K(\nu)}DK(\cc A)$$
where $wC^b(\cc A)$ stands for the category of quasi-isomorphisms
in $C^b(\cc A)$ and $K(wC^b(\cc A))$ is its Waldhausen $K$-theory.
The map on the left is induced by the map taking an object
$A\in\cc A$ to the complex concentrated in the zeroth degree and
the map $\nu$ is induced by the quotient functor $C^b(\cc A)\to
D^b(\cc A)$.

\begin{question}[The first Maltsiniotis conjecture~\cite{M}]
The map $K(\rho):K(\cc A)\to DK(\cc A)$ is a homotopy equivalence.
\end{question}

Let $\cc A$ admit an exact fully faithful functor $i:\cc A\to\cc
B$ with $\cc B$ an abelian category such that for any map $f$ in
$\cc A$ with $i(f)$ an epimorphism in $\cc B$ the map $f$ is an
epimorphism. This is the case when weak idempotent object split in
$\cc A$ (see~\cite{T}). Then $K(\tau)$ is a homotopy equivalence
by the Gillet-Waldhausen theorem~\cite[1.11.7]{T}. In this case,
the comparison conjecture is reduced to showing that $K(\nu)$ is a
homotopy equivalence.

Let us consider the composed map of spaces
   $$K(f):K(\cc A)\xrightarrow{K(\rho)}DK(\cc A)\xrightarrow{K(\mu)}K(\bb D^b(\cc A))$$
in which the right arrow is a homotopy equivalence. It is shown
in~\cite{Gar} that $K_0(f)$ is an isomorphism and hence is so
$K_0(\rho):K_0(\cc A)\to DK_0(\cc A)$.

I personally doubt that the comparison conjecture is true. This is
caused by a recent observation of To\"en and Vezzosi~\cite{TV}:
the obvious generalization of this conjecture to all Waldhausen
categories can not be true for obvious functoriality reasons. It
is true for the Waldhausen $K$-theory of spaces, for example.

Now we want to formulate a sort of Localization Theorem for the
$DK$-theory. We think that the following ingredients would be the
most reasonable to do that.

(1) One should first find the relevant notions of a thick exact
subcategory $\cc A$ of an exact category $\cc U$ and a quotient
exact category $\cc U/\cc A$ satisfying the obvious universal
property in $(ExCats)$.

(2) If $\cc A\subseteq\cc U$ is thick then so is $S_n\cc
A\subseteq S_n\cc U$ for every $n$.

(3) If $\cc A$ is thick and idempotent complete in $\cc U$ then
the sequence of bounded derived categories
   $$D^b(\cc A)\to D^b(\cc U)\to D^b(\cc U/\cc A)$$
is an exact sequence of triangulated categories, i.e. $D^b(\cc A)$
is the full triangulated category of $D^b(\cc U)$ on objects zero
in $D^b(\cc U/\cc A)$ and $D^b(\cc U)/D^b(\cc A)=D^b(\cc U/\cc
A)$.

The desired notions are suggested by Schlichting
in~\cite{Sch-noch} (see also his
Disser\-ta\-tions\-schrift~\cite{Schrift}). The conditions (1)-(2)
are satisfied for $\cc A\subseteq\cc U$ whenever $\cc A$ is a
``left or right s-filtering subcategory" in $\cc U$ in the sense
of~\cite{Sch-noch} (for brevity thick) and if, moreover, $\cc A$
is idempotent complete, then (3) is also valid. For example, any
filtering subcategory in the sense of Karoubi~\cite{Kar} or
Pedersen-Weibel~\cite{PW} is thick. Notice that if all categories
considered are abelian, then any thick subcategory is Serre.

\begin{question}[Localization]
If $\cc A$ is a thick and idempotent complete subcategory of an
exact category $\cc U$ then the sequence of exact categories $\cc
A\to\cc U\to\cc U/\cc A$ induces a homotopy fibration of spaces
   $$DK(\cc A)\to DK(\cc U)\to DK(\cc U/\cc A).$$
\end{question}

Localization would follow if we showed that the quotient functor
$\cc U\to\cc U/\cc A$ induces a homotopy equivalence
   $$i.\bb S.\bb S.(\cc A\subset\cc U)\to i.\bb S.\bb S.(0\subset\cc U/\cc A).$$
Indeed, we would have then the following commutative diagram
   $$\xymatrix{i.\bb S.\cc A\ar[r] & i.\bb S.\cc U\ar[r]\ar[d]&i.\bb S.\bb S.(\cc A\subset\cc U)\ar[d]\\
               &i.\bb S.\cc U/\cc A\ar[r]&i.\bb S.\bb S.(0\subset\cc U/\cc A)}$$
in which the first horizontal line is a homotopy fibration by
Corollary~\ref{333}. The right arrow is a homotopy equivalence and
$i.\bb S.\cc U/\cc A\to i.\bb S.\bb S.(0\subset\cc U/\cc A)$ is a
homotopy equivalence, too (for example by appealing again to
Corollary~\ref{333}).

\section{Addendum}

In this section we give the definition of a left pointed
d\'erivateur and its $K$-theory. The theory of d\'erivateurs was
developed by Grothendieck in~\cite{G}. Very close to d\'erivateurs
objects (the so-called ``homotopy theories" and ``systems of
diagrams categories") have been studied by Heller~\cite{H} and
Franke~\cite{F}. Since this paper mostly deals with the
d\'erivateur given by the hyperfunctor
   $$I\longmapsto\bb D^b(\cc A^I)$$
where $\cc A$ is an exact category we will only discuss
d\'erivateurs and its $K$-theory space although the analogous
$K$-theory can also be defined for systems of diagram categories
(see~\cite{Gar} for details).

\subsection{The axioms}

For the notions of the 2-category and 2-functor we refer the
reader to~\cite{Mac}. In what follows we use the term ``poset" as
an abbreviation of ``finite partially ordered set". The 2-category
of the posets (respectively the finite categories without cycles)
we shall denote by $\ord$ (respectively by $\dirf$).

Let $\di$ be a full 2-subcategory of the 2-category $\Cat$ of
small categories that contains the 2-category $\ord$. We assume
that $\di$ satisfies the following conditions:

\begin{enumerate}
\item $\di$ is closed under finite sums and finite products; \item
for any functor $f:I\to J$ in $\di$ and for any object $y$ of $J$,
      the categories $f/y$ and $f\setminus y$ are in $\di$.
\end{enumerate}

We shall also refer to $\di$ as a {\it category of diagrams}.

A {\it pred\'erivateur of the domain $\di$\/} or just a {\it
pred\'erivateur\/} is a functor
   \begin{equation}\label{dia}
    \bb D:\di^{\op}\to\CAT
   \end{equation}
from $\di$ to the category $\CAT$ of categories (not necessarily
small) satisfying the Functoriality Axiom below. So to each
category $I$ in $\di$ there is associated a category $\bb D_I$,
and to each map $f:I\to J$ in $\di$ a functor $f^*=\bb D(f):\bb
D_J\to\bb D_I$.

\begin{funcax}{\rm The following conditions hold:

\begin{itemize}

\item[$\diamond$] to each natural transformation $\phi:f\to g$ a
natural transformation $\phi^*:f^*\to g^*$ is associated and the
maps $f\to f^*$ and $\phi\to\phi^*$ define a functor from $\homa
IJ$ to the category of functors from $\bb D_J$ to $\bb D_I$;

\item[$\diamond$] if
   $$\xymatrix{K\ar[r]^f &I\ar@/^/[r]^g \ar@/_/[r]_{g'} &J\ar[r]^h &L}$$
are morphisms and $\phi:g\to g'$ is a bimorphism, then
$f^*\circ\phi^*=(\phi\circ f)^*$ and $\phi^*\circ
h^*=(h\circ\phi)^*$.

\end{itemize}
}\end{funcax}

A {\it morphism\/} $F:\bb D\to\bb D'$ between two pred\'erivateurs
$\bb D$ and $\bb D'$ consists of the following data:
\begin{enumerate}
\item for any $I\in\di$, a functor $F:\bb D_I\to\bb D'_I$;

\item for any map $f:I\to J$ in $\di$, $f^*F=Ff^*$;

\item for any bimorphism $\phi:f\to g$ in $\di$,
$\phi^*F=F\phi^*$.
\end{enumerate}

A morphism $F:\bb D\to\bb D'$ is an {\it equivalence\/} if for any
$I\in\di$ the functor $F:\bb D_I\to\bb D'_I$ is an equivalence of
categories.

Given $I\in\di$ and $x\in I$, let $i_{x,I}:0\to I$ be the functor
sending 0 to $x$. For $A\in\bb D_I$ let $A_x=i_{x,I}^*A$. Let us
consider the following axioms listed below.

\begin{isomax}{\rm
A morphism $f:A\to B$ in $\bb D_I$ is an isomorphism if and only
if $f_x:A_x\to B_x$ is so for all $x\in I$. }\end{isomax}

\begin{disjax}{\rm
(a) If $I=I_1\coprod I_2$ is a disjoint union of its full
subcategories $I_1$ and $I_2$, then the inclusions
$i_{1;2}:I_{1;2}\to I$ define an equivalence of categories
   $$(i^*_1,i^*_2):\bb D_I\lra{\sim}\bb D_{I_1}\times\bb D_{I_2}.$$

(b) $\bb D_\emptyset$ is a trivial category (having one morphism
between any pair of objects). }\end{disjax}

\begin{kanax}{\rm
The left homotopy Kan extension axiom says that for any functor
$f:I\to J$, the functor $f^*:\bb D_J\to\bb D_I$ has a left adjoint
$f_!:\bb D_I\to\bb D_J$. }\end{kanax}

\begin{baseax}{\rm
Let $f:I\to J$ be a morphism in $\di$ and $x\in J$. Consider the
diagram in $\di$
   $$\xymatrix{\ar @{}[dr] |{\swarrow\alpha_x}
               f/x \ar[d]_{p} \ar[r]^{j_x} & I \ar[d]^{f} \\
               0 \ar[r]_{i_{x,J}}        & J}$$
with $j_x$ a natural map and $\alpha_x$ the natural bimorphism.
The $\alpha_x$ induces a bimorphism $\beta_x:p_!j_x^*\to
i_{x,I}^*f_!$. The left base change axiom requires $\beta_x$ to be
an isomorphism. }\end{baseax}

\begin{defs}{\rm
A pred\'erivateur is called a {\it left d\'erivateur\/} if the
Functoriality Axiom, the Isomorphism Axiom, the Disjoint Union
Axiom, the Left Homotopy Kan Extension Axiom, and the Left Base
Change Axiom are satisfied. }\end{defs}

Let $F:\bb D\to\bb D'$ be a morphism between two left
d\'erivateurs, and let $f:I\to J$ be a map in $\di$. Consider the
adjunction maps $\alpha:1\lra{}f^*f_!$ and $\beta:f_!f^*\lra{}1$.
Denote by $\gamma_{F,f}$ the composed map
   $$f_!F\xrightarrow{f_!F\alpha}f_!Ff^*f_!=f_!f^*Ff_!\xrightarrow{\beta Ff_!}Ff_!.$$
$F$ is {\it right exact\/} if $\gamma_{F,f}$ is an isomorphism and
the following compatibility relations hold:
   $$F\alpha_{\bb D}=f^*(\gamma_{F,f})\circ\alpha_{\bb D'}F
    \textrm{ and }
    F\beta_{\bb D}=\beta_{\bb D'}F\circ\gamma^{-1}_{F,f}f^*.$$

The d\'erivateurs we work with are also to meet some extra
conditions. A map $f:I\to J$ in $\di$ is a {\it closed (open)
immersion\/} if it is a fully faithful inclusion such that for any
$x\in J$ the relation $\Hom(I,x)\ne\emptyset$
($\Hom(x,I)\ne\emptyset$) implies $x\in I$.

\begin{defs}{\rm
A left d\'erivateur is {\it pointed\/} if the following conditions
hold:

(1) for any closed immersion $f:I\to J$ in $\di$, the structure
functor $f_!$ possesses a left adjoint $f^?$;

(2) for any open immersion $f:I\to J$ in $\di$, the structure
functor $f^*$ possesses a right adjoint $f_*$;

(3) for any open immersion $f:I\to J$ in $\di$ and any object
$x\in J$, the base change morphism of the diagram
   $$\xymatrix{\ar @{}[dr] |{\nearrow\gamma_x}
               f\setminus x \ar[d]_{q} \ar[r]^{l_x} & I \ar[d]^{f} \\
               0 \ar[r]_{i_{x,J}}        & J}$$
yields an isomorphism $\delta_x:i^*_{x,I}f_*\to q_*l_x^*$.
}\end{defs}

If $\bb D$ is a left pointed d\'erivateur, then $\bb D_I$ has a
zero object for any $I\in\di$.

\subsection{The $S.$-construction and $K$-theory space}

Throughout this section $\bb D$ is assumed to be a left pointed
d\'erivateur (of the domain $\di$).

Let $\Box\in\di$ be the poset $\Delta^1\times\Delta^1$ and let
$\la\subset\Box$ be the subposet $\Box\setminus (1,1)$. Let
$i_{{}_\ulcorner}:\la\to\Box$ be the inclusion. An object $A$ of
$\bb D_\Box$ is called {\it cocartesian\/} if the canonical
morphism $i_{{}_{\ulcorner!}}i_{{}_\ulcorner}^*A\to A$ is an
isomorphism.

Given $0\le i<j<k\le n$ let
   \begin{equation}\label{aijk}
    a_{i,j,k}:\Box\to\Ar\Delta^n
   \end{equation}
denote the functor defined as:
   $$(0,0)\longmapsto(i,j),\ \ \ (0,1)\longmapsto(i,k),\ \ \ (1,0)\longmapsto(j,j),
     \ \ \ (1,1)\longmapsto(j,k).$$
For any integer $n\geq 0$, denote by $S_n\bb D$ the full
subcategory of $\bb D_{\Ar\Delta^n}$ that consists of the objects
$X$ satisfying the following two conditions:
\begin{itemize}
\item[$\diamond$] for any $i\le n$, $X_{(i,i)}$ is isomorphic to
zero in $\bb D_0$;

\item[$\diamond$] for any $0\le i<j<k\le n$, $a_{i,j,k}^*X$ is a
cocartesian square if $n>1$.
\end{itemize}
The definition of $S_n\bb D$ is similar to that of $S_n\cc C$,
where $\cc C$ is a Waldhausen category~\cite{W}. $S_0\bb D$ is the
full subcategory of zero objects in $\bb D_0$. The category
$S_1\bb D$ consists of the objects $X\in\bb D_{\Delta^2}$ with
$X_0$ and $X_2$ isomorphic to zero.

For any object $I$ of $\di$, we denote by $\bb D(I)$ the left
pointed d\'erivateur defined as $\bb D(I)_J=\bb D_{I\times J}$.
Let $\bb S_n\bb D_I=S_n\bb D(I)$. Then $\bb S_n\bb D$ is a left
pointed d\'erivateur (see~\cite{Gar}). There results a simplicial
left pointed d\'erivateur
   $$\bb S.\bb D:\Delta^n\longmapsto \bb S_n\bb D.$$

Consider the following simplicial category:
   $$S.\bb D:\Delta^n\longmapsto S_n\bb D.$$
For any $n\geq 0$, let $iS_n\bb D$ denote the subcategory of
$S_n\bb D$ whose objects are those of $S_n\bb D$ and whose
morphisms are isomorphisms in $S_n\bb D$, and let $i.S_n\bb D$ be
the nerve of $iS_n\bb D$. We obtain then the following
bisimplicial object:
   $$i.S.:\Delta^m\times\Delta^n\longmapsto i_mS_n\bb D.$$

\begin{defs}{\rm
The {\it Algebraic $K$-theory\/} for a small left pointed
d\'erivateur $\bb D$ of the domain $\di$ is given by the pointed
space (a fixed zero object 0 of $\bb D_0$ is taken as a basepoint)
   $$K(\bb D)=\Omega|i.S.\bb D|.$$
The $K$-groups of $\bb D$ are the homotopy groups of $K(\bb D)$
   $$K_*(\bb D)=\pi_*(\Omega|i.S.\bb D|)=\pi_{*+1}(|i.S.\bb D|).$$
}\end{defs}

Denote by $(\textit{Left\ pointed\ d\'erivateurs})$ the category
of left pointed d\'erivateurs and right exact functors. Then the
map
   $$(\textit{Left\ pointed\ d\'erivateurs})\lra{K}(\textit{Spaces})$$
is functorial.

Let $\bb D$ be a left pointed d\'erivateur. Denote by $\bb E_0$
the full subcategory in $\bb D_\Box$ consisting of the cocartesian
squares $E\in\bb D_\Box$ with $E_{(1,0)}$ isomorphic to zero. If
we replace $\bb D$ by $\bb D(I)$, we define the category $\bb E_I$
similar to $\bb E_0$. One obtains a left pointed d\'erivateur $\bb
E$. It is equivalent (by a right exact map) to the d\'erivateur
$\bb S_2\bb D$~\cite[6.2]{Gar}.

\end{document}